\documentclass[proc]{edpsmath}
\usepackage{fullpage}
\usepackage[utf8]{inputenc}
\usepackage[T1]{fontenc}
\usepackage{amsmath,amsfonts,amssymb,mathtools}
\usepackage{stmaryrd}
\usepackage{graphicx, array}
\usepackage{mathrsfs,dsfont}

\usepackage{color}

\usepackage{enumerate}

\usepackage{hyperref}

\newtheorem{theo}{Theorem}[section]
\newtheorem{propo}[theo]{Proposition}
\newcommand{\E}{\mathbb{E}}
\newcommand{\R}{\mathbb{R}}
\newcommand{\PP}{\mathbb{P}}

\newcommand{\N}{\mathbb{N}}
\newcommand{\dt}{\Delta t}

 \newcommand{\saut}[1]{\vspace{#1\baselineskip}}
\newcommand{\abs}[1]{\left|#1\right|}
\newcommand{\set}[1]{\left\{#1\right\}}
\newcommand{\pare}[1]{\left ( #1\right )}

\newcommand{\espace}[1]{\mathds{#1}}

\begin{document}

\title{Analysis and simulation of rare events for SPDE\small s}\thanks{This work was supported by the CEMRACS 2013, CIRM, Luminy}

\author{Charles-Edouard Br\'ehier}\address{Universit\'e Paris-Est, CERMICS (ENPC),  6-8-10 Avenue Blaise Pascal, Cit\'e Descartes ,  F-77455 Marne-la-Vall\'ee, France; \email{brehierc@cermics.enpc.fr}}\secondaddress{INRIA Paris-Rocquencourt, Domaine de Voluceau - Rocquencourt, B.P. 105 - 78153 Le Chesnay, France;\\ \email{mathias.rousset@inria.fr}}
\author{Maxime Gazeau}\address{INRIA Lille - Nord Europe, 
Parc Scientifique de la Haute-Borne, Park Plaza b\^atiment A, 40 avenue Halley, 59650 Villeneuve d'Ascq Cedex, France; \email{maxime.gazeau@inria.fr}}
\author{Ludovic Gouden\`ege}\address{F\'ed\'eration de Math\'ematiques de l'\'Ecole Centrale Paris, CNRS, Grande voie des vignes, 92295 Ch\^atenay-Malabry, France; \email{goudenege@math.cnrs.fr}}
\author{Mathias Rousset}\sameaddress{1,2}

\begin{abstract}
In this work, we consider the numerical estimation of the probability for a stochastic process to hit a set $B$ before
reaching another set $A$. This event is assumed to be rare. We consider \emph{reactive trajectories} of the stochastic Allen-Cahn partial differential evolution equation (with double well potential) in dimension $1$. Reactive trajectories are defined as the probability distribution of the trajectories of a stochastic process, conditioned by the event of hitting $B$ before $A$. We investigate the use of the so-called Adaptive Multilevel Splitting algorithm in order to estimate the rare event and simulate reactive trajectories. This algorithm uses a \emph{reaction coordinate} (a real valued function of state space defining level sets), and is based on (i) the selection, among several replicas of the system having hit $A$ before $B$, of those with maximal reaction coordinate; (ii) iteration of the latter step. We choose for the reaction coordinate the average magnetization, and for $B$ the minimum of the well opposite to the initial condition. We discuss the context, prove that the algorithm has a sense in the usual functional setting, and numerically test the method (estimation of rare event, and transition state sampling).
\end{abstract}

\maketitle



\section{Introduction}

\subsection{The goal}
This paper focuses on the metastable states and associated rare events of reversible stochastic gradient systems in general dimension (finite or infinite). Reversible stochastic gradient systems can be used as a paradigmatic model to study physical metastability, that it is to say systems exhibiting very stable states (in terms of time scale), which are however very unlikely. In $\R^d$ with the usual euclidean structure, typical reversible stochastic gradient systems are given by diffusions solutions of Stochastic Differential Equations (SDEs) of the form:
\begin{equation}\label{eq:sde}
dX_t=-\nabla V(X_{t})dt+\sqrt{2\epsilon}\, d B_t,
\end{equation}
where $t \mapsto B_t \in \R^d$ is a standard Wiener process (Brownian Motion). If $V:\R^d \to \R $ is a smooth map such that:
\begin{enumerate}[(i)]
\item the reversible invariant Gibbs measure
  \begin{equation}
    \label{eq:gibbs}
    \mu_{\epsilon}(dx) := \frac{1}{Z_\epsilon} {\rm e}^{-\frac{1}{\epsilon} V(x) } \, d x, 
  \end{equation}
is a probability measure (the equilibrium distribution) for the appropriate normalisation constant $Z_\epsilon$;
\item the local minima of $V$ form a countable set of isolated points;
\end{enumerate}
then it is well-known that when $\epsilon \to 0$, the equilibrium distribution concentrates on the global minimum $$ V(x_0) = \inf_{x \in \R^d} V(x),$$ while all other local minima become metastable states (precise statements are known as the Freidlin-Wentzell large deviation theory, see below). A typical physical context leading to models of the form~\eqref{eq:sde} is given by stochastically perturbed molecular dynamics. In the latter, $V$ is the classical interaction potential energy between the atoms of molecules, and $X_t$ represents the position of the atoms after a mass weighted linear change of coordinates. Moreover, the stochastic perturbation (thermostat) is assumed to be Markovian, with very strong viscosity (the so-called ''overdamped'' limit), and adimensional temperature proportional to $\epsilon$. Finally, we assume we have quotiented out the continuous symmetries of the system so that $V$ has indeed only isolated local minima.

\saut{0.25}

Let us now define two closed subset $A \subset \R^d$ and $B \subset \R^d$; for instance, one may keep in mind the case where $X_0=x_0$, $ A = \bigcup_{n \geq 2} B(x_n,\delta) $ and $B = B(x_1,\delta)$, the latter being small balls centered at the local minima $x_n, n \geq 1$. We next define the hitting time of $A$ (or equivalently $B$) as
$$\tau_A \equiv \tau_A ^\epsilon =\inf\left\{t\geq 0; X_t \in A\right\}.$$ We are interested in computing the rare event probability
\begin{equation}
  \label{eq:proba}
 0 <  \mathbb{P}(\{ \tau_B < \tau _A\}) \ll 1,
\end{equation}
as well as sampling according to the \emph{reactive path} ensemble, defined as the conditional probability distribution on trajectories 
\begin{equation}
  \label{eq:react_traj}
  {\rm Law} ( X_t, \, t \geq 0 \vert \{ \tau_B < \tau _A\} ).
\end{equation}

\subsection{Infinite dimension}
This work will also focus on gradient systems in \emph{infinite dimension}. As a consequence the stochastic differential equation~\eqref{eq:sde} will be replaced by a Stochastic Partial Differential Equation perturbed with a space-time white noise - see also the abstract formulation \eqref{SPDE2}:
\begin{equation}\label{SPDE1}
\begin{cases}
\dfrac{\partial x(t,\lambda)}{\partial t}=\gamma \dfrac{\partial^2 x(t,\lambda)}{\partial \lambda^2}-\nabla V(x(t,\lambda))+\sqrt{2\epsilon}\dfrac{\partial^2\omega(t,\lambda)}{\partial t\partial \lambda}, \quad t>0 \text{ and } \lambda \in (0,1)\\[0.25cm]
\dfrac{\partial x(t,\lambda)}{\partial \lambda}|_{\lambda=1}=\dfrac{\partial x(t,\lambda)}{\partial \lambda}|_{\lambda=0}=0\\[0.25cm]
x(0,\lambda)=x_0(\lambda),
\end{cases}
\end{equation}
where $x_0$ is a given initial condition. The second line expresses the homogeneous Neumann boundary conditions at $\lambda=0$ and $\lambda=1$. The deterministic forcing term $\nabla V$ is a sufficiently regular vector field on a Banach space $H$, and $\mathcal{V}(x) = \int_{0}^1V(x(\lambda))d \lambda$ is a non-quadratic potential energy. 
The problem of sampling according to the rare event~\eqref{eq:proba} and the simulation of reactive trajectories~\eqref{eq:react_traj} remains unchanged.

\subsection{Splitting algorithms}
In finite dimension, more efficient algorithms than the plain Monte-Carlo simulation of independent realizations of~\eqref{eq:sde} were developed introducing a real valued function on the state space, usually called \emph{reaction coordinate}, and denoted by
\[
\xi: \R^d \to \R.
\]
 We will assume that $X_0=x_0$ with 
\[
\sup_{x \in A} \xi(x) \leq z_A < \xi(x_0) < z_B \leq \inf_{x \in B} \xi(x).
\]
The general key point of splitting methods to compute the rare event $\{ \tau_B < \tau _A\}$ consists in simulating $n_{rep}$ replicas of~\eqref{eq:sde}, and to duplicate with higher probability the trajectories with highest maximal level $\xi$. It turns out that this can be done in a consistent way. In Section \ref{AMS}, a general adaptive algorithm inspired by \cite{CerGuy07, CerGuyLel11} and an associated unbiased estimator of the probability of $\{ \tau_B < \tau _A\}$ will be presented, and numerically tested in Section \ref{num}.

One objective of this paper is to generalize this approach to an infinite dimensional setting.

\subsection{Large deviations and Transition State Theory}

The behavior of the solutions to either \eqref{eq:sde} or \eqref{SPDE1} in the limit~$\epsilon \to 0$ is well explained by the theory of large deviations \cite{FreiWen}. On a finite time window $[0,T]$, it can then be checked that on the space of trajectories (with uniform topology) 
\begin{align}\label{diff:pot}
W(x,y) = \inf_{T \geq 0 }\inf_{x_0=x,x_T=Y} I(x) = \sup_{x_0=x,x_\infty=y }V(x_t) - V(x),
\end{align}
with \emph{good rate function}
\[
I(x) := \frac{1}{4}\int_{0}^T \abs{\nabla V (x_t) - \dot{x}_t}^2 \, dt,
\]
and where extrema are taken over smooth trajectories. The left-hand side of \eqref{diff:pot} is the so-called \emph{quasi-potential function} defined for any $(x,y)\in \R^d$ as the minimal cost of forcing the system to go from $x$ to $y$ in an indeterminate time. For gradient systems, the latter formula shows that it is the lowest energy barrier that needs to be overcome in order to reach $y$ from $x$. The quasi-potential yields a rationale to compute rare events related to \emph{exit times}, since for any $\delta >0$, we have the general formula
\[
\lim_{\epsilon \to 0} \mathbb{P}_{X_0=x} ( {\rm e}^{( \inf_{y \in \partial A } W(x,y) - \delta )/\epsilon} < \tau^\epsilon_A <  {\rm e}^{( \inf_{y \in \partial A } W(x,y) + \delta )/\epsilon}  ) = 1.
\]
This shows that if $\inf_{y \in \partial A } W(x,y) < \inf_{y \in \partial B } W(x,y)$, then the event $\{ \tau_B < \tau _A\}$ is indeed a rare event in the limit $\epsilon \to 0$. Tools of Potential Theory are used to go further allowing to compute the so-called prefactor; see \cite{Bovier} for the general methodology and results in finite dimension. In the double well situation, where $x_0$ denotes saddle point and $x_{-/+}$ the minima, Kramer's law holds in the sense that
\[
\espace{E}_{x_-}(\tau_{B_{\delta}(x_+)}) = \frac{2\pi}{\abs{\lambda_1(x_0)}}\sqrt{\frac{\abs{\det(\nabla^2 V(x_0))}}{\det(\nabla^2 V(x_-))}}e^{(V(x_0) -V(x_-))/\epsilon}\left[1+ \mathcal{O}\left(\epsilon^{1/2}\abs{\log \epsilon}^{3/2}\right)\right]
\]
where $\tau_{B_{\delta}(x)}$ is a small ball of radius $\delta$ centered at $x$.

\saut{0.25}

In the infinite dimensional setting, the choice of functional spaces is essential. In \cite{CheMi}, the large deviation result is given in H\"older-like spaces $\mathcal{C}^{2\beta,0}([0,1])$ for $\beta\in]0,1/4[$, with Dirichlet boundary conditions. In~\cite{Barret, BerGen}, the Kramer's law is given respectively with respect to Sobolev and H\"older norms. But this kind of information is asymptotic, merely theoretical, and the direct computation of the above quantities is not possible in general. This justifies the introduction of efficient numerical methods.

\section{Finite and infinite dimensional models}
In this Section, we introduce two different mathematical models with metastability, which are linked together through
their energies. The first model is a coupled system of $N$ stochastic differential equations called 
the overdamped Langevin equation. It models the evolution of $N$ interacting particles driven by random forces. The second
one is a stochastic partial differential equation known as the Allen-Cahn equation. The main motivation of this study is 
that the second equation can be viewed, under an appropriate scaling, as the limit of the SDE when $N$ goes to infinity.

\saut{0.5}
 
We first introduce known results on these two problems and necessary notations especially on the functional spaces. 
We also present various discretization methods that can be implemented to obtain numerical approximations of these processes.

\subsection{The finite dimensional model}
Throughout this paper, we will use the following terminology: an \textit{atom} denotes a particle and a system of $N$ atoms is defined as a \textit{molecule}, with $ N \in \N$ and $N\geq 2$. We model a molecule with $N$ atoms moving on a line and submitted to three kinds of forces:
\begin{enumerate}
\item each atom $i \in \llbracket 1,N \rrbracket$ is confined by a potential $V :\R \to \R$ with at least two wells. One particular 
interesting case, considered in this paper, is the symmetric double-well potential $V(x)=\frac{x^4}{4}- \frac{x^2}{2}$;
\item each atom $i$ interacts with its two nearest neighbors $i-1$ and $i+1$ through a quadratic potential;
\item the movement of each atom is perturbed by a small Gaussian white noise in time. The processes acting on different atoms are mutually independent.
\end{enumerate}
 The noise perturbation is represented by a $N$-dimensional Wiener process $W^N=(W_1^N,\ldots, W_N^N)$ whose covariance is given 
for any $u,v\in \mathbb{R}^N$ by
\begin{equation}
\E\bigl[\left\langle W^N(t),u\right\rangle\left\langle W^N(s),v\right\rangle\bigr]=\min(t,s)\left\langle u, v\right\rangle,
\end{equation}
for any $s,t\in\R^+$ and where  $\left\langle . , . \right\rangle$ denotes the canonical scalar product on $\mathbb{R}^N$. 
The position of the atom $i$ at time $t$ is denoted by $X_i^N(t)$. The configuration of the molecule $(X^N_i)_{1\leq i\leq N}$ satisfies the following system of Stochastic Differential Equations for $t> 0$ and any $i \in \left\{1, \cdots, N\right\}$
\begin{align}\label{SystSDE}
dX_i^N(t)&=\gamma N(X_{i-1}^N(t)+X_{i+1}^N(t)-2X^N_i(t))dt\\&-\frac{1}{N}V'(X_i^N(t))dt+\sqrt{2\epsilon}dW_i^N(t) \nonumber,
\end{align}
with homogeneous discrete Neumann type boundary conditions $X_{N+1}^N:= X_{N}^N$, $X_{0}^N:= X_{1}^N$ and for an initial configuration $(X_i(0))_{1\leq i\leq N}=(x_i^N)_{1\leq i\leq N}$. 

\saut{0.25}

The positive parameter $\gamma>0$ determines the rigidity between the atoms.
The value of $\gamma$ is of importance in the study of metastability. Indeed the nature of some stationary points of the energy defined below is modified, as well as the reactive trajectories between the metastable states.
The parameter $\epsilon>0$ 
represents the level of the noise acting on the molecule. Physically it can be interpreted as the temperature. One of
our concerns is to investigate numerically the evolution of the molecule when $\epsilon$ decreases to $0$
and we aim to propose efficient algorithms to estimate the transitions rates from one stable state to another. The total
energy $\mathcal{E}_\gamma^N$ of the system is given for any $x^N=(x_1^N,\ldots,x_N^N)\in\mathbb{R}^N$
\begin{equation}\label{defEnergyN}
\mathcal{E}_\gamma^N(x^N)=\frac{\gamma}{2N}\sum_{i=0}^{N}N^2(x_{i+1}^{N}-x_{i}^{N})^2+\frac{1}{N}\sum_{i=1}^{N}V(x_i^N),
\end{equation}
with the boundary conditions $x_{N+1}^N:= x_{N}^N$, $x_{0}^N:= x_{1}^N$.
The system \eqref{SystSDE} can be rewritten as the \textit{overdamped Langevin dynamics} in $\mathbb{R}^N$ associated with
the energy $\mathcal{E}_\gamma^N$. 
Indeed, \eqref{SystSDE} can be rewritten as
\begin{equation}\label{overdampedN}
dX^N(t)=-\nabla \mathcal{E}_\gamma^N(X^N(t))dt+\sqrt{2\epsilon}dW^N(t).
\end{equation}
The regularity of the potential $V:\mathbb{R}\rightarrow \mathbb{R}$ is important for the global well-posedness of this system.
We will give sufficient regularity conditions for it.
\begin{rmrk}
This choice of boundary terms is linked to the homogeneous Neumann boundary conditions imposed on the SPDE below. 
Moreover the scaling of the different contributions in the energy with respect to the size $N$ in the molecule is related to the convergence when $N\rightarrow +\infty$.
\end{rmrk}

\subsection{Infinite dimensional model}
Let $x:[0,1]\rightarrow \mathbb{R}$ be a function of class $\mathcal{C}^1$ on the space interval $[0,1]$. For any positive integer $N$, we denote $\Delta x=\frac{1}{N+1}$ and $x_i^N=x(i\Delta x)$. The generic space variable in the interval $[0,1]$ is denoted by $\lambda$. Then 
\begin{align*}
\mathcal{E}_{\gamma}^{N}(x^N_1,\ldots,x^N_{N})\xrightarrow[N \to +\infty]{}
\mathcal{E}_{\gamma}^{\infty}(x)=\frac{\gamma}{2}\int_{0}^{1}|x'(\lambda)|^2d\lambda+\int_{0}^{1}V(x(\lambda))d\lambda
\end{align*}
as Riemann sums. Inspired by this convergence for smooth functions, we define a new energy on $H^1(0,1)$, the space of square integrable functions defined on $(0,1)$ having a derivative (in the distributional sense) belonging to $L^2(0,1)$:
\begin{equation}\label{defEnergyInf}
\mathcal{E}_{\gamma}^{\infty}(x)=\frac{\gamma}{2}\int_{0}^{1}|x'(\lambda)|^2d\lambda+\int_{0}^{1}V(x(\lambda))d\lambda.
\end{equation} 
We recall that a Sobolev embedding ensures that $H^1(0,1)$ is continuously embedded into $\mathcal{C}([0,1])$; thus the function $V\circ x$ is well-defined and continuous on $[0,1]$ for any $x\in H^1(0,1)$.

\saut{0.25}

The stochastic partial differential equation linked to this energy is the Allen-Cahn equation - see \eqref{SPDE1}. In an abstract form, it is written
\begin{equation}\label{SPDE2}
  \begin{cases}
    dX(t)=AX(t)dt - \nabla V(X(t))dt+\sqrt{2\epsilon}dW(t),\\
    X(0)=x_0.
  \end{cases}
\end{equation}
The interpretation of \eqref{SPDE2}, as a stochastic evolution equation on (an appropriate subset of) the Hilbert space $H=L^2(0,1)$, follows the classical framework of \cite{[DPZ]}. The operator $A$ denotes the realization on $H$ of the Laplace operator $\frac{\partial^2}{\partial \lambda ^2}$ with homogeneous Neumann boundary conditions. Accordingly its domain is
\[
\mathcal{D}(A) = \left\{ \phi \in H^2(0,1), \phi'(0)= \phi'(1)=0\right\}.
\]
The non-linear mapping $F(x)(\lambda)=-\nabla V(x(\lambda))$,  for almost all $\lambda\in(0,1)$, is defined on the domain $D(F)=L^p(0,1)$, where $p$ is such that $F(x)\in L^2(0,1)$ for any $x\in L^p(0,1)$. In the case of the Allen-Cahn potential, $V(x)=\frac{x^{4}}{4}-\frac{x^2}{2}$, and one can take $p=6$. More general assumptions for $V$ are possible:
\begin{enumerate}
	\item $V$ is of class $C^3$;
	\item $V$ is convex at infinity, i.e. there exists two constants $C>0$ and $c > 0$ such that for any $\abs{x}>c $ then
	\[
	V''(x) > C > 0;
	\]
	\item There exists two constants $C>0$ and $\alpha \geqslant 2$ such that for any $x \in \R$ and $j=0,1,2,3$
	\[
	\abs{V^j(x)} \leqslant C (\abs{x}^{2\alpha-1} +1).
	\]
\end{enumerate}
The noise perturbation in \eqref{SPDE2} is induced by a centered, Gaussian space-time white noise. Its covariance 
satisfies for any times $s,t$ and any positions $\lambda,\mu$
\[\E[\omega(t,\lambda)\omega(s,\mu)]=\min(t,s)\times \min(\lambda,\mu).\]
The noise $W$ is given by a cylindrical Wiener process. Formally, it is defined as the following series
\[W(t)=\sum_{k\in\mathbb{N}}W_k(t)e_k,\]
where $(e_k)_{k\in\mathbb{N}}$ is any complete orthonormal system of the Hilbert space $H$, and $(W_k)_{k\in\mathbb{N}}$ is a sequence of independent standard scalar Wiener processes. The convergence of this series does not hold in $H$, but only on larger spaces in which $H$ is embedded thanks to a linear Hilbert-Schmidt operator. 
The eigenvalues of $A$ are then $\lambda_n = (n\pi)^2$ for $n\in\N$
and the corresponding eigenfunctions are $e_n(\lambda)=\sqrt{2}\cos(n\pi \lambda)$ for $n=1,2 \cdots$, and $e_0(\lambda)=1$; they satisfy
\[
\sum_{k=1}^{+\infty} \frac{1}{\lambda_k^{\alpha}} < \infty \Leftrightarrow \alpha >1/2.
\]
If we denote $\widetilde{e}_k = \frac{e_k}{\lambda^{1/2+\delta}_k}$ for $k\geq 1$ and $\widetilde{e}_0=e_0$, then the series
\[\widetilde{W}(t)=\sum_{k\in\mathbb{N}}W_k(t)\widetilde{e}_k,\]
converges. The white noise is thus well defined in $H^{-1/2 - \delta}$
for any $\delta >0$.
The heat kernel will play the role of the Hilbert Schmidt operator and give a rigorous meaning to the stochastic convolution
\[
W_{A}(t) = \int_0^t e^{(t-s)A}dW(s)
\]
in $H^{1/2 -\delta}$, for any $\delta>0$.

\begin{prpstn}
Let $\delta >0$, and fix a suitable parameter $p$. For any final time $T>0$ and continuous initial data $x_0\in L^p(0,1)$ satisfying the Neumann boundary condition, then Equation \eqref{SPDE2} admits a unique 
mild solution in $L^2(\Omega, C( [0,T] \times [0, 1]) )$:
\begin{equation}\label{mild}
X(t)=e^{tA}x_0+\int_{0}^{t}e^{(t-s)A}F(X(s))ds+\sqrt{2\epsilon}\int_{0}^{t}e^{(t-s)A}dW(s),
\end{equation}
such that for every $t\in[0,T]$ we have $X(t)\in L^p(0,1)$. Moreover it admits a modification which
is of class $C^{1/4-\delta}$ in time and $C^{1/2-\delta}$ in space. The solution of the linearized
equation is in $H^{1/2-\delta}$ for any $\delta >0$. Moreover, we have the following stronger estimate:
$$\E[\sup_{t\in[0,T],\lambda\in[0,1]}|x(t,\lambda)|^2]\leq C(T)<+\infty.$$
\end{prpstn}
This proposition is proved using classical truncation and approximation arguments \cite{Barret}.
Moreover, model \eqref{SystSDE} can be seen as an approximation with finite differences of Equation \eqref{SPDE2}.
For the convergence of $X^N$ to $X$, see for instance \cite{Barret}.

\subsection{Invariant distribution}
In the case of the double-well potential, some critical points of the deterministic dynamics ($\epsilon=0$) in both finite and infinite dimensions are easily identified.

For the overdamped Langevin dynamics in dimension $N<+\infty$, two global minima are the molecules
$$x_N^{\pm}=\pm(1,\ldots,1).$$ 
Another critical point is $$x_N^0=(0,\ldots,0);$$
however as explained in Section \ref{Bif} the precise nature of $x_N^0$, given by the signature of the Hessian $D^2\mathcal{E}_{\gamma}^{N}(x^0)$, depends on $\gamma$ and $N$. Moreover when $\gamma$ changes new critical points may appear.

In the infinite dimensional case, the corresponding global minima are the constant functions
$$x_\infty^{\pm}(\lambda)=\pm 1, \quad \text{for all }\lambda\in[0,1].$$
Remark that these functions satisfy the homogeneous Neumann boundary conditions.

We also have a critical point such that $x_\infty^0(\lambda)=0$. 

\saut{0.25}

A distinction must be done to analyze the invariant laws according to the dimension. There is always ergodicity and a formula of Gibbs type, but with respect to a different type of reference measure which is the Lebesgue measure in $\R^N$ and a suitable (degenerate) Gaussian measure in $H$ in infinite dimension.

\saut{0.25}

\textit{Finite dimension: } The dynamics is given by the overdamped Langevin equation \eqref{overdampedN}, with the energy $\mathcal{E}_{\gamma}^{N}$. Assumptions on $V$ ensure that the probability measure defined with
\[
\mu_{\gamma}^{N,\epsilon}(dx) = \frac{1}{Z_{\gamma}^{N,\epsilon}}\exp\left(-\frac{\mathcal{E}_{\gamma}^{N}(x)}{\epsilon}\right)dx,
\]
is the unique invariant law, where $dx$ denotes Lebesgue measure on $\R^N$ and $Z_{\gamma}^{N,\epsilon}$ is a normalization constant. Moreover, the stationary process is ergodic and reversible with respect to this invariant measure. There is another way to write $\mu_{\gamma}^{N,\epsilon}$, when $\gamma>0$, as a Gibbs measure with respect to the measure
\[
\nu_{\gamma}^{N,\epsilon}(dx)=\exp\left(-\frac{\gamma\mathcal{D}^{N}(x)}{\epsilon}\right)dx,
\]
using the following decomposition of the energy into kinetic and potential energy of the atoms in the molecule:
\begin{equation}\label{eq:energy}
\begin{gathered}
\mathcal{E}_{\gamma}^{N}=\gamma\mathcal{D}^{N}+\mathcal{V}^{N}(x)\\
\mathcal{D}^{N}(x)=\frac{1}{2N}\sum_{i=1}^{N}N^2(x_{i+1}^{N}-x_{i}^{N})^2\\
\mathcal{V}^{N}(x)=\frac{1}{N}\sum_{i=1}^{N}V(x_i^N).
\end{gathered}
\end{equation}
The measure $\nu_{\gamma}^{N,\epsilon}$ cannot be normalized as a probability measure. It is a degenerate Gaussian measure, which can be seen as the unique (up to a multiplicative constant) invariant measure of the linearized process - i.e. with $V=0$. Using a suitable change of variables, it is expressed as the product of the one-dimensional Lebesgue measure and of a non-degenerate Gaussian probability measure. The degeneracy is explained by the choice of Neumann boundary conditions: from \eqref{SystSDE} with $V=0$ (linearization) and $\epsilon=0$ (deterministic case) we see that the average $\xi(X_1,\ldots,X_N)=\frac{1}{N}\sum{i=1}^{N}X_i$ is preserved. When $\epsilon>0$, there is no dissipation on this mode, which is solution of $d\xi_t=\sqrt{2\epsilon/N}dB^N(t)$ where $B^N(t)=\frac{1}{N}\sum_{i=1}^{N}W_i^N(t)$.

\saut{0.5}

\textit{Infinite dimension: } The above discussion needs to be adapted. Recall that if $X(t)$ is solution of the SPDE then almost surely, the quantity $\mathcal{E}_\gamma(X(t))$ is not finite for $t>0$. However, there is an expression of the invariant law, with a density with respect to a reference measure $\nu_{\gamma}^{\infty,\epsilon}$. The linear SPDE - with $V=0$ - is written
\[
dY(t)=\gamma AY(t)+\sqrt{2\epsilon}dW(t).
\]
We decompose the Hilbert space $H$ into $H=H_0+H_0^{orth}$ with $H_0=\text{Span}(e_0)$. Then if we denote $Y(t)=Y_0(t)+Y_{0}^{orth}$, the above linear SPDE is decomposed into two decoupled equations
\begin{gather*}
dY_0(t)=\sqrt{2\epsilon}dW_0(t),\\
dY_{0}^{orth}(t)=\gamma AY_{0}^{orth}(t)+\sqrt{2\epsilon}dW_{0}^{orth}(t).
\end{gather*}
$H_0^{orth}=\text{Span}\left\{e_i;i\geq 1\right\}$ is a stable subspace of $A$, and $A$ is invertible on $H_0^{orth}$, with inverse denoted $A_{orth}^{_1}$. The equation on $Y_0$ admits Lebesgue measure $\nu_0(dx_0)=dx_0$ as unique (up to a multiplicative constant) invariant measure. The equation on $Y_{0}^{orth}$ admits a unique invariant law, which is a Gaussian measure on $H_{0}^{orth}$, denoted by $\nu_{0}^{orth}(dx_{0}^{orth})$; it is centered and its covariance operator is $\frac{1}{\gamma\epsilon}A_{orth}^{_1}$. Therefore the linearized SPDE admits the following invariant measure defined in $H$:
\[
\nu_{\gamma}^{\infty,\epsilon}(dx)=\nu_0(dx_0) \times \nu_{0}^{orth}(dx_{0}^{orth}).
\]
It is worth noting that for any $p\in[2,+\infty[$ the Banach space is included in the support of $\nu_{\gamma}^{\infty,\epsilon}$.
Assumptions on the potential $V$ now ensure that the SPDE \eqref{SPDE2} admits a unique invariant probability measure when $\gamma>0$, denoted by $\mu_{\gamma}^{\infty,\epsilon}$, with the expression
\[
\mu_{\gamma}^{\infty,\epsilon}(dx) = \frac{1}{Z_{\gamma}^{\infty,\epsilon}} \exp\left\{-\frac{1}{\epsilon} \int_0^1 V(x(\lambda))d\lambda\right\} \nu_{\gamma}^{\infty,\epsilon}(dx),
\]
where $Z_{\gamma}^{\infty,\epsilon}$ is a normalization constant.

\begin{rmrk}
The case of homogeneous Dirichlet boundary conditions is studied in \cite{Otto:13, Weber:10}. There is an interpretation of the Gaussian measure $\nu_{\gamma}^{\infty,\epsilon}$ as the law of the Brownian Bridge, with a renormalization with respect to the parameters $\gamma$ and $\epsilon$.
\end{rmrk}

Transitions between the stable equilibrium points appear when the temperature parameter $\epsilon$ is positive; the typical time can be expressed thanks to the Kramers law.

\subsection{Discretization of the stochastic processes}
The aim of this part is to introduce somehow classical schemes to solve numerically the Langevin overdamped equation and 
the stochastic Allen-Cahn equation. For the finite dimensional system we consider the well known Euler scheme
 \[
 X^{n+1} = X^n -  \nabla\mathcal{E}^N_{\gamma} (X^n) \dt+ \sqrt{2 \epsilon \dt} G_n 
 \]
where $G_n$ are standard Gaussian random variables. This scheme is proved to be
of order $1/2$ in $L^p(\Omega)$ norms for any $p\geq 2$. Regarding the numerical approximation of the stochastic Allen-Cahn
equation, various numerical schemes may be proposed that are based on deterministic scheme. 
We consider here finite difference schemes but other methods may be used efficiently such as 
finite elements. Spectral methods seem to be a little harder 
to implement. Indeed one has to use Chebychev polynomials instead of the fast Fourier transform
because Neumann conditions are imposed on the boundary.
The nonlinear term is handled using a splitting method. 
The basic idea of splitting methods (see \cite{Strang:68} is to
approach the exact flow $L(t)X_0$ of the nonlinear equation by means of a truncation of 
the Baker-Campbell-Hausdorff formula. Let us denote by $S(t)X_0$ the solution 
of the stochastic heat equation
\begin{align*}
\begin{cases}
dX(t)&= \gamma \partial^2_xX(t)dt + \sqrt{2\epsilon}dW(t)\\
X(0) &= X_0,
\end{cases}
 \end{align*}
and $T(t)Y_0$ the solution of the Bernoulli differential equation 
\begin{align*}
\begin{cases}
\partial_tY(t)&= - \left(Y^3(t) -Y(t)\right) \\
Y(0) &= Y_0.
\end{cases}
 \end{align*}
The first equation can be easily solved numerically using
a second-order semi-implicit scheme to insure the  unconditional stability of the scheme.
The second equation is exactly solvable and its solution is given by
\[
Y(t) = \dfrac{Y_0}{\sqrt{Y_0^{2} + (1-Y_0^{2})\exp(-2t)}}.
\]
The Lie method consists to approximate the exact solution of the stochastic Allen-Cahn equation
$L(t)X_0$ by either one of the two methods $S(t)T(t)X_0$ or $T(t)S(t)X_0$ which correspond to 
a composition of the two previous flows. When $\epsilon =0$, 
it is well known that the order of convergence of this scheme is $1$ in time. Adding the white noise, 
the strong order of convergence usually drop to $1/4$. The numerical Lie scheme reads as follows:
 \begin{align} \left\{\begin{array}{ll}X_{n+1/2} - X_n = \dfrac{\gamma\dt}{4} \partial_x^2 (X_{n+1/2}+ X_n)  + \dfrac{1}{2}\sqrt{2 \epsilon \dt} G^n \\[0.3cm]
	  X_{n+1}   =  \dfrac{X_{n+1/2}}{\sqrt{X_{n+1/2}^{2} + (1-X_{n+1/2}^{2})\exp(-2\dt)}} 
		\end{array}\right. \end{align}
where $G^n =\sum_{j=0}^J G_j^n e_j $ and $G_j^n\sim\mathcal{N}(0,1)$ are iid. $J$ is a truncation parameter for the expansion of the noise in a basis of the Hilbert space adapted to the linear operator.

Higher order schemes may be constructed considering more terms in the Baker-Campbell-Hausdorff formula.
For time independent and deterministic operators, a method to construct even order symplectic integrator is proposed in \cite{Yoshida:90}. 
A natural extension of such schemes is:
\begin{equation}\label{eq::scheme}
\begin{aligned} \left\{\begin{array}{ll}X_{n+1/3} - X_n = \dfrac{\gamma\dt}{4} \partial_x^2 (X_{n+1/3}+ X_n)  + \dfrac{1}{2}\sqrt{2 \epsilon \dt} G^n \\[0.3cm]
	  X_{n+2/3}   =  \dfrac{X_{n+1/3}}{\sqrt{X_{n+1/3}^{2} + (1-X_{n+1/3}^{2})\exp(-2\dt)}} \\[0.5cm]
	 X_{n+1} - X_{n+2/3} = \dfrac{\gamma\dt}{4} \partial_x^2 (X_{n+1} + X_{n+2/3})  + \dfrac{1}{2}\sqrt{2 \epsilon \dt } G^n\end{array}\right. \end{aligned}
\end{equation}
The equation is then discretized with finite differences.
We have implemented this method for our simulations; however the analysis of its order of convergence is an open question.

\section{Bifurcations}\label{Bif}
The total energy $\mathcal{E}^N_{\gamma}$ depends on a parameter $\gamma$ that determines 
the influence of the gradient energy with respect to the potential one. It is refereed 
as a bifurcation parameter and a small change in its value may suddenly change the behavior
of the dynamical system. It can be seen studying the nature and the number of the critical points of the energy.
In both finite and infinite dimension, one may thus expect that the reactive trajectories do not experiment the same paths according to the value of $\gamma$.
To illustrate this, let us consider the energy functional for $N=2$ as defined in~\eqref{defEnergyN} and study its critical points. The approach here is inspired by \cite{Bouchet:09, Corvellec:12}.
\begin{propo}
For $N=2$, the local maximum $0$ is unstable in the sense that it degenerates into 
a local maximum at the bifurcation parameter $\gamma=1/8$ in the direction given by the vector $(1,-1)$.
Then for $\gamma < 1/8$, two saddle points appear at $\pm \sqrt{1-8\gamma}$.
A new change of regime occurs at $\gamma =1/12$ and these two saddle points degenerate to local minima; four new saddle points appear when $\gamma<1/12$.

\saut{0.5}

For $N=4$, the saddle point $0$ is also unstable. It degenerates at the bifurcation parameters 
\[
\gamma = \left\{\frac{1}{16(2-\sqrt{2})}, \frac{1}{16}, \frac{1}{16(2+\sqrt{2})}\right\}.
\]
A bifurcation occurs in the direction $(1,0,0, 1)$.
\end{propo}
For $N=2$, bifurcations are studied by two approaches: a direct computation and introducing suitable normal forms.
Normal forms are simplified functionals exhibiting the same structure of critical points than the full problem and describing the phase transitions in a neighbourhood of a bifurcation value. This is done for $N=2$ to emphasize the importance of normal forms in cases when explicit computations can not be performed (for example $N=4$).

\saut{0.25}

\emph{Direct computation:} We denote $\gamma  = \gamma/4$. The Jacobian matrix is
\[
\nabla\mathcal{E}^2_{\gamma} (x,y) = \left(\begin{array}{c}
x^{3}-x + \gamma (x-y)\\
y^{3}-y - \gamma (x-y)
\end{array}\right).
\]
Summing the two equations, we easily obtain that the critical points satisfy the system
\begin{align*}
\left\{ \begin{array}{ll}
(x+y)(x^{2}-xy+y^{2}-1) = 0 \\[0.2cm]
  x^{3}-x + \gamma (x-y)=0.
	\end{array}\right.
\end{align*}
Let us now determine the real roots of this system with respect to $\gamma$:
 \begin{enumerate}
	 \item If $x=-y$ then $x$ satisfies
\begin{align}
x(x^{2} + (2\gamma-1))=0. \label{E1}
\end{align}
Consequently $x=0$ or $x=\pm \sqrt{(1-2\gamma)}$.
\item If  $x \neq -y$ then 
\begin{align}\label{E2}
\left\{ \begin{array}{ll}
(x^{2}-xy+y^{2}-1) = 0  \\[0.2cm]
x^{3} - x+\gamma (x-y)=0.
	\end{array}\right.
\end{align}
From the first equation of system \eqref{E2}, $y = (x \pm \sqrt{4-3x^{2}})/2$. Plugging this expression into the second equation of \eqref{E2} leads to an equation for $x$
\begin{align*}
(x-1)(x+1)(x^4+x^2+ \gamma^2)=0,
\end{align*}
whose roots are given by $\pm 1$ and
\[
\alpha(\gamma) =\pm\sqrt{\frac{1-\gamma\pm\sqrt{(\gamma+1)(1-3\gamma)}}{2}}.
\]
 \end{enumerate}
We now sum up the critical points of the energy according to some values of the  bifurcation parameter $\gamma$
\[
\gamma\geqslant \frac{1}{2}\Longrightarrow\left\{\begin{array}{c}
x=y=0\\
x=y=\pm1
\end{array}\right.
\]
and 
\[
\frac{1}{3}\leqslant\gamma\leqslant \frac{1}{2}\Longrightarrow\left\{\begin{array}{c}
 x=y=0\\
x=y=\pm1\\
x=-y=\pm\sqrt{1-2\gamma}\\
\end{array}\right.
\]
If $\gamma \leqslant \frac{1}{3}$ then $\alpha(\gamma)$ is real valued and the critical points are
\[
\left\{\begin{array}{c}
x=y=0\\
x=y=\pm1\\
x=-y=\pm\sqrt{1-2\gamma}\\
x=\alpha(\gamma) \text{ and } y = \frac{x(x^{2}-1+\gamma)}{\gamma}.
\end{array}\right.
\]
Let us now determine their nature and see how a small change of $\gamma$ may influence the variations of the energy.
The Hessian of the energy $\mathcal{E}^2_{\gamma}$ is given by
\[
\mathcal{H}\mathcal{E}^2_{\gamma} (x,y) = \left(\begin{array}{cc}
3x^{2}-1 + \gamma & -\gamma\\
-\gamma&3y^{2}-1 + \gamma 
\end{array}\right).
\]
At the point $(0,0)$, its spectrum is given by $Sp = \{-1, 2\gamma -1\}$. When $\gamma > 1/2$, the point $(0,0)$ is a saddle point while it degenerates into a local maximum for $\gamma \leqslant 1/2$. Another change of regimes occurs at the points
$(\pm \sqrt{1-2\gamma}, \mp \sqrt{1-2\gamma})$. Indeed
\[
\mathcal{H}\mathcal{G} (\pm\sqrt{1-2\gamma},\mp\sqrt{1-2\gamma}) = \left(\begin{array}{cc}
2-5\gamma & -\gamma\\
-\gamma&2 -5\gamma 
\end{array}\right) \Longrightarrow Sp =  \{2-4\gamma, 2-6\gamma\}.
\]
Therefore for $1/2 \geqslant \gamma \geqslant 1/3$, these two points are saddles while they degenerate to local minima when $\gamma < 1/3$.

\saut{0.25}

\begin{figure}[h]
\centering
\includegraphics[width=12cm, height =10cm]{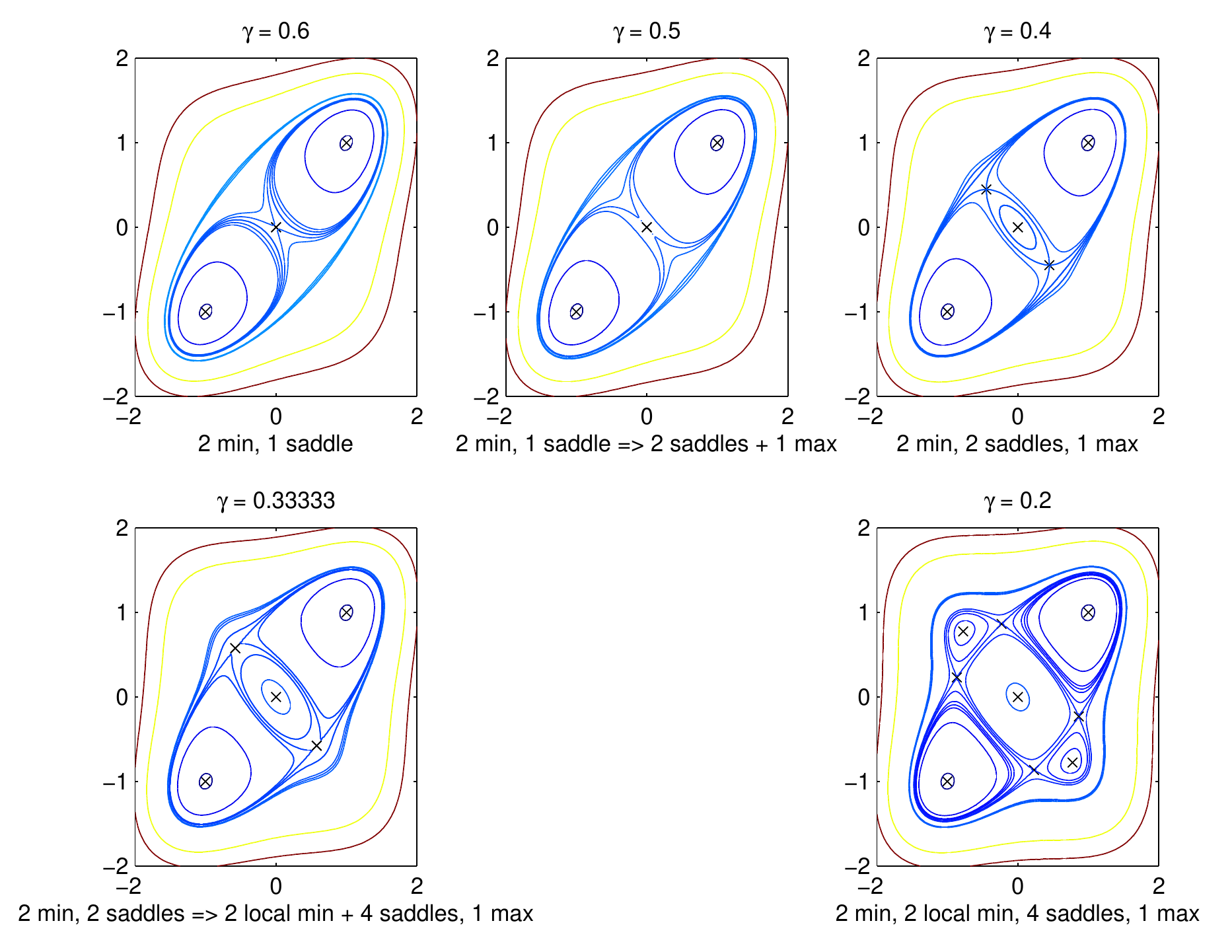}
\caption{Bifurcations for different values of $\gamma$. Plot of $\mathcal{E}^2_{\gamma}$.}
\label{Fig:Bif}
\end{figure}

\emph{Normal form:}
Bifurcations at $(0,0)$ are expected to be described by a normal form in a direction orthogonal to the eigenvector $(1,1)$. In this simple case, the orthogonal space is easily identified and is $(1,-1)$. We decompose $q$ along this basis
\[
q = A\begin{pmatrix}
1 \\ 
1
\end{pmatrix} + \rho\begin{pmatrix}
1 \\ 
-1
\end{pmatrix}.
\]
Plugging $q$ into the energy, one easily gets 
$\mathcal{E}^2_{\gamma}(\rho,A) = \mathcal{G}_{0} + \mathcal{G}_{1}$,
where
\begin{align*}
\left\{ \begin{array}{ll}
\mathcal{G}_{0}(A) =  A^4/2 -A^2  \\[0.2cm]
\mathcal{G}_{1}(\rho,A) = \rho^4/2-  \rho^2(1- 2\gamma) +3\rho^2A^2 .
\end{array}\right.
\end{align*}
We introduce a new quantity 
\begin{align}
G_1(A) : = \min_{\rho} \mathcal{G}_{1}(\rho,A)  = \begin{cases}
-\frac{1}{2}+2\gamma+3A^2-2\gamma^2-6\gamma A^2-\frac{9}{2}A^4  &\mbox{if } \abs{A} < \sqrt{\frac{1 - 2\gamma}{3}}\\
0 &\mbox{otherwise }
\end{cases}
\end{align}
and the normal form is given by
\begin{align*}
G(A) :=\mathcal{G}_{0}(A) + G_1(A) = \begin{cases} -1/2+2\gamma+2A^2-2\gamma^2-6A^2\gamma-4A^4  &\mbox{if } \abs{A} < \sqrt{\frac{1 - 2\gamma}{3}}\\
A^4/2 -A^2  &\mbox{otherwise}.
\end{cases}
\end{align*}
Then the nature of the normal form changes at $\gamma=1/2$. For $\gamma >1/2$, the normal form is locally concave around 
zero which is a local maximum. For $\gamma <1/2$, the critical points of this functional are now $0$ and $\pm  \sqrt{1-3\gamma}/2$. For a small perturbation $\gamma_{\epsilon} = 1/2 - \epsilon, \ 0 <\epsilon \ll 1$, 
the minimum of $G$ is reached at $0$ and
\[
G(0) = -\frac{1}{2}+2\gamma_{\epsilon}-2\gamma_{\epsilon}^2.
\]
Then we conclude that $-G(A)$ describes the bifurcation of $\mathcal{E}_{\gamma}^2$ around $(0,0)$ in the direction $(1,-1)$. 
Fig \ref{Fig:Bif} plots the isolines 
of the energy $\mathcal{E}^2_{\gamma}$ for different values of $\gamma$. 

For $N=4$ and $x= (x_1, x_2, x_3, x_4)$, the Jacobian matrix of $\mathcal{E}^4_{\gamma} (x)$ is
\[
\nabla\mathcal{E}^4_{\gamma} (x) = \begin{pmatrix}
\frac{1}{4}\left(x_1^{3}-x_1\right) - 4\gamma (x_2-x_1) \\[0.2cm]
\frac{1}{4}\left(x_2^{3}-x_2\right) - 4\gamma (x_3-2x_2+x_1)   \\[0.2cm]
\frac{1}{4}\left(x_3^{3}-x_3\right) - 4\gamma (x_4-2x_3+x_2) \\[0.2cm]
\frac{1}{4}\left(x_4^{3}-x_4\right) + 4\gamma (x_4-x_3)
\end{pmatrix}
\]
and its Hessian 
\[
\mathcal{H}\mathcal{E}^4_{\gamma} (x) = \begin{pmatrix} \frac{3}{4}x_1^2-\frac{1}{4}+4\gamma & -4\gamma& 0 & 0\\
  -4\gamma &  \frac{3}{4}x_2^2- \frac{1}{4}+8\gamma& -4\gamma& 0 \\ 
	0 & -4\gamma &  \frac{3}{4}x_3^2- \frac{1}{4}+8\gamma & -4\gamma \\ 0 & 0& -4\gamma &  \frac{3}{4}x_4^2- \frac{1}{4}+4\gamma
\end{pmatrix}
\]
Obvious critical points are $0, \pm 1$. Direct computations are more complex than in two dimensions and we make 
use of the normal forms to study bifurcations. Let us first note that the spectrum of $\mathcal{H}\mathcal{E}^4_{\gamma} (0)$ is 
\[Sp=\left\{-\frac{1}{4}, -\frac{1}{4}+8\gamma, (8+4\sqrt{2})\gamma-\frac{1}{4}, (8-4\sqrt{2})\gamma-\frac{1}{4}\right\}.\]
The associated eigenvectors are respectively
\begin{align*}
\begin{pmatrix}
1\\
1\\
1\\
1
\end{pmatrix}, \quad \begin{pmatrix}
1\\
-1\\
-1\\
1
\end{pmatrix}, \quad \begin{pmatrix}
1\\
-1-\sqrt{2}\\
1+\sqrt{2}\\
-1
\end{pmatrix},\quad \begin{pmatrix}
1\\
-1+\sqrt{2}\\
1-\sqrt{2}\\
-1
\end{pmatrix}  
\end{align*}
It is obvious that the saddle point $0$ degenerates into different types of saddle points at  
\[\gamma = \left\{\frac{1}{16(2-\sqrt{2})}, \frac{1}{16}\right\} \]
and finally to a local maximum at $\gamma = 1/(16(2+\sqrt{2}))$. 
A full bifurcation diagram for $N=4$ can be found in \cite{Berglund:06} but for a different energy. In our case, bifurcations
do not appear for the same values of $\gamma$ but this diagram gives a good insight of what may happen 
in higher dimension than $2$.
\begin{figure}[h]
\centering
\includegraphics[width=12cm, height =8cm]{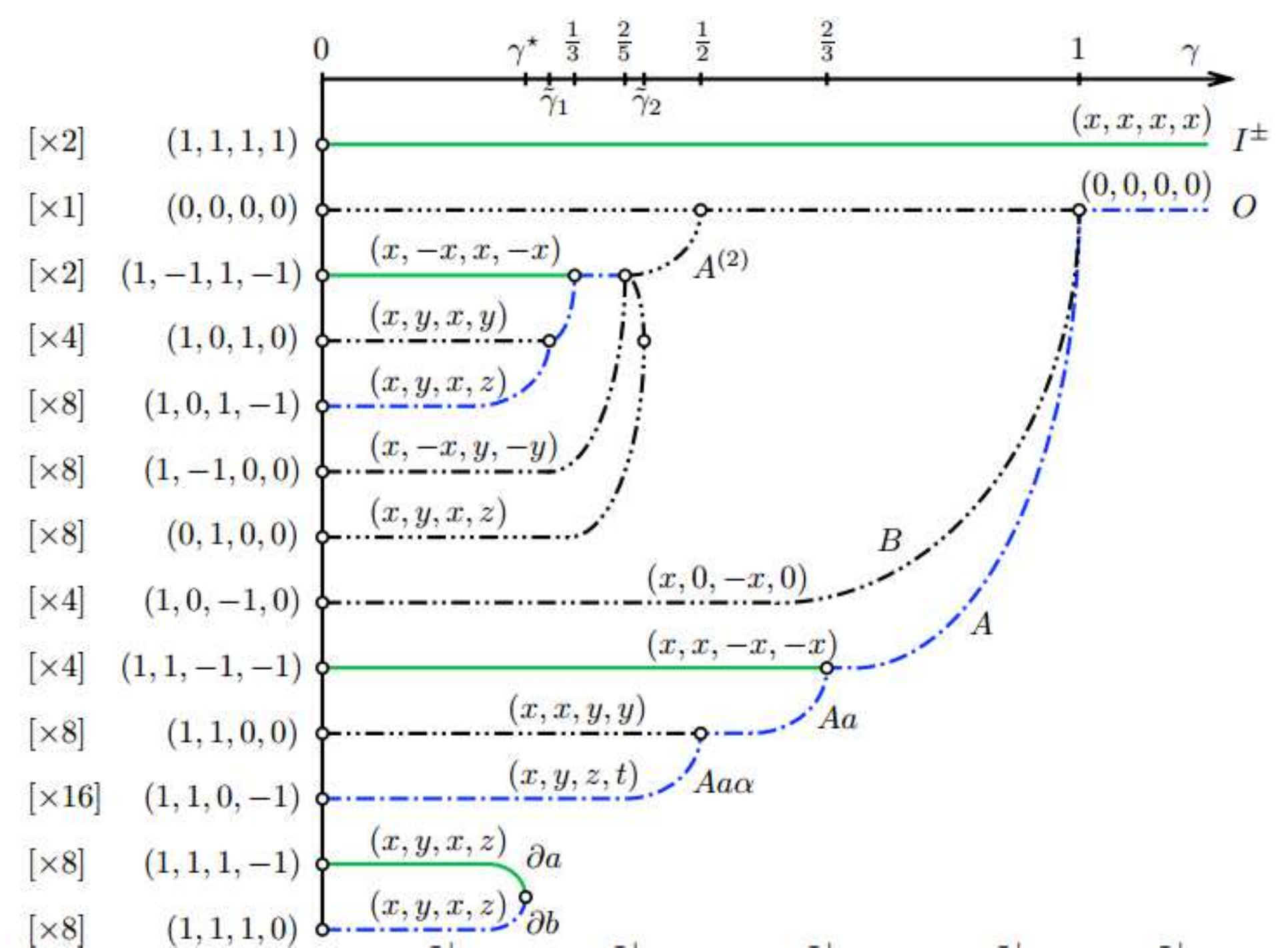}
\caption{From \cite{Berglund:06}. Bifurcations for different values of $\gamma$.}
\label{Bif:4atomes}
\end{figure}
We study bifurcations in the orthogonal of the eigenvector corresponding to the biggest eigenvalues. Accordingly we consider $x = A e_1 + \rho'$, 
where $e_1 = (1, -1-\sqrt{2}, 1+\sqrt{2}, -1)$ and $\rho'$ belongs to the orthogonal of $e_1$. Obviously this situation is more complex 
than in dimension $2$ since the orthogonal space is now of dimension three. Let us first consider the
case where $\rho' = \rho(1,0,0,1)^t$. Then,
\[
\mathcal{E}^4_{\gamma}(x) = \frac{9}{4}A^4+ \frac{1}{8}\rho^4-A^2-\frac{1}{4}\rho^2+48\gamma A^2+\frac{3}{4}A^2\rho^2+32\gamma A^2\sqrt{2}+4\gamma\rho^2+\frac{3}{2}A^4\sqrt{2}-\frac{1}{2}A^2\sqrt{2}.
\]
Denoting $\mathcal{G}_1(\rho, A) = \frac{1}{8}\rho^4-\frac{1}{4}\rho^2+\frac{3}{4}A^2\rho^2+4\gamma\rho^2$, the critical points are given by $0$ and $\pm\sqrt{-16\gamma-3A^2+1}$. Thus  
\begin{align*}
G_1(A) := \min_{\rho}\mathcal{G}_1(\rho, A) = \begin{cases}
-\frac{1}{8}+\frac{3}{4}A^2+4\gamma-\frac{9}{8}A^4-12\gamma A^2-32\gamma^2 &\mbox{if } \abs{A} < \sqrt{\frac{1 -16\gamma }{3}}\\
0 &\mbox{otherwise}.
\end{cases}
\end{align*}
Therefore the normal form is given by
\begin{align*}
G(A) = \begin{cases}
\frac{9}{8}A^4-\frac{1}{8}-\frac{1}{4}A^2+4\gamma+36\gamma A^2-32\gamma^2+32\gamma A^2\sqrt{2}+\frac{3}{2}A^4\sqrt{2}-\frac{1}{2}A^2\sqrt{2}  &\mbox{if } \abs{A} < \sqrt{\frac{1 -16\gamma }{3}}\\
\frac{9}{4}A^4-A^2+48\gamma A^2+32\gamma A^2\sqrt{2}+\frac{3}{2}A^4\sqrt{2}-\frac{1}{2}A^2\sqrt{2}  &\mbox{otherwise.}
\end{cases}
\end{align*}
If $\gamma>1/16$, then $G'(A) = 9A^3-2A+96\gamma A+64\gamma A\sqrt{2}+6A^3\sqrt{2}-A\sqrt{2}$ and the only critical point is zero which is a local minimum. At $\gamma= 1/16$ a bifurcation occurs and 
\[
G'(A) = \frac{9}{2}A^3-\frac{1}{2}A+72\gamma A+64\gamma A\sqrt{2}+6A^3\sqrt{2}-A\sqrt{2}.
\]
There are now three critical points: zero is a local maximum while the two other critical points
are local minima. Therefore in the direction $e_1$, a local maximum degenerates into a minimum. Thanks to the use of normal forms we are able to partly describe bifurcations. 
Much analysis has to be done in other directions considering linear combinations of orthogonal vectors.

\saut{0.25}

In the case $N=+\infty$, such a discussion is possible, in order to study the nature of the critical point $x_\infty^0$.

\section{Adaptive Multilevel Splitting}\label{AMS}

The goal of splitting methods is to simulate $n$ replicas (copies) of reactive trajectories. Let us recall that reactive trajectories are defined by the conditional distribution \eqref{eq:react_traj}, the distribution of the SDE or SPDE dynamics~\eqref{eq:sde}-\eqref{SPDE2} conditioned by the event $\{ \tau_B < \tau _A\} $.

Loosely speaking, the small probability of the rare event as in enforced using a birth-death mechanism as follows:
\begin{itemize}
\item Launch multiple replicas subject to the reference dynamics. 
\item Kill the replicas having the smallest maximal level, as given by the reaction coordinate (fitness) continuous mapping $\xi: \mathbb{R}^{N} \text{ or } H\rightarrow \mathbb{R}$.
\item Replicate the other replicas.
\end{itemize}
 
In the literature, the following cases have been studied:
\begin{itemize}
\item $N=1$, this is the classical AMS algorithm proposed in \cite{CerGuy07};
\item $N>1$, multiple replica algorithm based on AMS, studied in \cite{CerGuyLel11}.
\end{itemize}

The infinite dimensional case $N=\infty$ seems to be treated for the first time in the present work.

We will use the following \emph{reaction coordinate} (fitness mapping) defined as the \emph{average position or magnetization} $\xi: L^2(0,1) \to \R$:
\[
 \xi(x) := \int_{0}^{1}x(\lambda) \, d \lambda \, \in \R .
\]

For the finite dimensional model, a discrete version of the magnetization is used, based on a quadrature formula for the integral. For instance,
\[
 \xi^N(x^N) :=\frac{1}{N}\sum_{i=1}^{N}x_i^N \, \in \R .
\]

We will present the algorithm in \emph{discrete time}, with time index $k \in \mathbb{N}$, so that it can directly be applied to the numerical discretization.

Consider now $n_{\rm rep}$ replicas of the system dynamics. We thus consider the \emph{(stopping) hitting time} of \emph{open levels} of $\xi$:
\[
 t_{z}(X) := {\rm inf} \set{ k \geq 0 \, \vert \,  \xi(X_k) > z }  .
\]
The \emph{rare event} of interest is then rewritten as
\[
\set{ t_{z_B}(X) \leq  \tau_A(X) },
\]
and we wish to compute its associated \emph{small probability}. Recall that $\tau_A(X)={\rm inf} \set{ k \geq 0 \, \vert \, X_k\in A }$. If we take $A={x;  \xi(x)<z_A }$, we see that the $\tau_A(X)={\rm inf} \set{ k \geq 0 \, \vert \,  \xi(X_k) < z_A }$.

Typically for the Allen-Cahn equation, $\xi(X_{0})\in(z_A,z_B)$, $z_A = -1 + \delta$, $z_B =1-\delta$, with some $\delta > 0$. Since we have
$$\xi(x_{\infty}^{-})=-1, \quad  \xi(x_{\infty}^{0})=0, \quad \xi(x_{\infty}^{+})=+1,$$
the magnetization gives some useful piece of information along the transition from one metastable state to the other.

In principle the algorithme can be generalized to \emph{any continuous} $\xi: E \to \R$ on a Polish state space $E$.

First, let us fix
\[
\mathcal{T} \subset \mathbb{N},
\]
the set of time indices for which the replica states will be kept in memory. We then denote by
\[
k \mapsto X_k^{(n,q)}, \quad n = 1 \dots n_{\rm rep}
\] 
indexed by $q\geq 1$ the iteration index of the algorithm, which is different from the \emph{time index} $k\geq 0$.  

For $q=1$ (initial condition), $(k \mapsto X_k^{(n,1)})_{1 \leq n \leq N }$ are iid and stopped at 
\[
{\rm min}\pare{ \tau_A(X^{(n,1)} ),t_{z_B} (X^{(n,1)})} ,
\]
i.e. when the magnetization is either below level $z_A$ or above $z_B$.

Then iterate on $q\geq 1$ as follows:
\begin{enumerate}
\item $q$ being given, compute:
\[
\begin{cases}
\displaystyle {\rm Max}^{(n,q)} := \mathop{ \rm max}_{k \geq 0} \xi(X^{(n,q)}_k) \\
  \displaystyle N^{\rm killed}_q := \mathop{\rm argmin}_{n}  {\rm Max}^{(n,q)}  \in [1, N],
\end{cases}
\]
the replica $N^{\rm killed}_q$ (we assume that it is unique) with minimal maximal (''min-max'') level ${\rm Max}^{(N^{\rm killed}_q,q)}$.
\item Fix $m \in \mathbb{N}_\ast$ some . Choose $N^{\rm new}_q$ uniformly in $\set{ 1, \dots, N} - \set{N^{\rm killed}_q}$ and consider the time (for a small $\delta > 0$):
  \begin{equation*}
    \tau^{(q)} := \inf \left\{  k \in \mathcal{T}  \vert t_{ {\rm Max}^{(N^{\rm killed}_q,q)}}( X^{(N^{\rm new}_q,q)}  ) \geq k \right\},
  \end{equation*}
the first time when the branching replica $N^{\rm new}_q$ has reached the maximum level ${\rm Max}^{(N^{\rm killed}_q,q)}$ of all the killed replicas.
\item Kill the information of replica $N^{\rm killed}_q$. Copy the path of replica $N^{\rm new}_q$ until $\tau^{(q)} $ and then re-sample the remaining path with Markov dynamics until either level $z_A$ or $z_B$ is reached.
\item Stop at $q=Q_{\rm iter}$ all replicas have reached level $z_B$.
\end{enumerate}

Then the general principle of the AMS algorithm can be stated as follows. Assume $n_{\rm rep} \to +\infty $. Then:
\begin{enumerate}
\item A path of a replica at the end of the algorithm is distributed according to the ''reactive trajectory''~\eqref{eq:react_traj}:
\[
{\rm Law} \pare{X_{t}, \, t \geq 0 \, \big \vert \,  t_{z_B}(X)<\tau_A(X) }
\]
\item The quantity
 \[
(1 - 1/n_{\rm rep})^{Q_{\rm iter}}
\]
is a convergent estimator of $\PP( t_{z_B}(X)<\tau_A(X))$.
\end{enumerate}

In a work in preparation \cite{Nous}, we are in fact able to prove the following unbiased property, for a Markov chain in any Polish state space $E$:
\begin{prpstn}
 Let $(X_k)_{k \geq 0}$ and $(X_{\ast,k})_{k \geq 0}$ two i.i.d. copies of the Markov chain. Assume that for any initial condition $x \in E$
\[
\PP_x \left ( \mathop{ \rm max}_{k \geq 0} \xi(X_{k \wedge \tau_A(X)})  = \mathop{ \rm max}_{k \geq 0} \xi(X_{\ast,k \wedge \tau_A(X_\ast)}) \right ) =0.
\]
Then we have the unbiased estimation
\[
\E \left ( (1 - 1/n_{\rm rep})^{Q_{\rm iter}} \right ) = \PP( t_{z_B}(X)<\tau_A(X)).
\]
\end{prpstn}
There also exists a variant of the above algorithm where $k_rep>1$ replicas are killed at each iteration. A suitably defined estimator also satisfies the unbiased property.

\section{Numerical results}\label{num}
We now want to show the performance of the algorithm in two directions: the estimation of the transition probability and the approximation of the reactive trajectories when $\gamma$ varies.

\subsection{Estimation of the probability}
With the notations of the previous section, we estimate the probability $\PP( \tau_B(X) \leq \tau_A(X))$ where the 
sets $A$ and $B$ correspond to the levels $z_A=-0.99$, $z_B=0.99$. The Allen-Cahn equation is discretized thanks to the scheme \eqref{eq::scheme} with a time-step $\Delta t=0.01$, and a finite difference discretization with a regular grid with mesh size $\Delta x=0.02$. The parameter $\gamma$ is equal to $1$  and we consider the initial condition $X_0(\lambda)=-0.8$.

We perform $N_{MC}=100$ realizations of the algorithm to compute a Monte-Carlo approximation of the expectation of the estimator. More precisely, we consider different choices of $n_{rep}\in\left\{50,100,200,1000\right\}$ and $\epsilon=0.05$, so that the probability of the transition is approximately $0.005$. In Table \ref{tab::1} below, we give an estimated probability, and an empirical standard deviation computed with the Monte-Carlo approximation.

Confidence intervals are obtained when this standard deviation is divided by the square root of the number of realizations: in Table \ref{tab::1} below the empirical standard deviation should then be divided by $10$.

\begin{table}[h]
\begin{tabular}{|c|c|c|}
  \hline
  $n_{rep}$ & estimated probability & empirical standard deviation \\
  \hline
  50 & 0.00516 & 0.00163 \\
  100 & 0.00514 & 0.00128 \\
  200 & 0.00502 & 0.000774\\
  1000 & 0.00501 & 0.000350\\
  \hline
\end{tabular}
\label{tab::1}
\caption{Estimated probability and empirical standard deviation obtained via the AMS algorithm.}
\end{table}

We observed that the precision is improved when we $n_{rep}$ increases. 
Another useful comparison is given in Table \ref{tab::2}, where we compare the results for $(n_{rep}=100,N_{MC}=1000)$ and $(n_{rep}=1000,N_{MC}=100)$. 

\begin{table}[h]
\begin{tabular}{|c|c|c|c|c|}
\hline
$n_{rep}$ & $N_{MC}$ & estimated probability & empirical standard deviation & time for one realization \\
\hline
100 & 1000 & 0.00502 & 0.00166 & 67 s\\
1000 & 100 & 0.00501 & 0.000350 & 676 s \\
\hline
\end{tabular}
\label{tab::2}
\caption{Computational time and variance for two values of $n_{rep}$.}
\end{table}

In both cases we obtain the same precision, with approximately the same required computational time. Two arguments are then in favor of choosing the smallest $n_{rep}$ in this situation: first, the estimator is unbiased for every value of $n_{rep}$, so that taking $n_{rep}$ very large is not necessary; second, it is easy to save computational time with a parallelization of the Monte-Carlo procedure. Parallelization inside the AMS algorithm could also help, and further research is necessary in this direction.

Finally, we compare the performance of the AMS algorithm with a direct Monte-Carlo procedure: we run independent trajectories solving the Allen-Cahn equation, and count $1$ if $\tau_B<\tau_A$, $0$ if $\tau_A<\tau_B$, and average over the realizations. The computation of one trajectory only takes about $0.4$ s, and if we run $10^6$ independent replicas, the Monte-Carlo procedure gives an estimated probability $0.00507$, with an empirical standard deviation $7.12 10^{-5}$. Compared to the result with $n_{rep}=200$ in Table \ref{tab::1}, we see that for the same precision the AMS algorithm is between $3$ and $4$ times faster than a direct Monte-Carlo method, for the approximation of a probability of order $5.10^{-3}$. Moreover, the smaller the probability becomes, the better the AMS algorithm should be.



Now we study the dependence of the transition probability with respect to the temperature parameter $\epsilon$. For each value of $\epsilon$, we choose the same discretization parameters i.e. $\Delta t=0.01, \Delta x=0.02$ and $\gamma=2$. Moreover, we take $n_{rep}=100$ as well as $N_{MC}=100$ realizations in order to compute an empirical probability with a Monte-Carlo procedure. Results are given with the empirical standard deviation of the estimator in Table \ref{tab:3}

\begin{table}[h]
\begin{tabular}{|c|c|c|}
\hline
$\epsilon$  & estimated probability & empirical standard deviation\\
\hline
0.30 & 0.0831 & 1.33 $10^{-2}$\\
0.10 & 0.0276 & 5.11 $10^{-3}$\\
0.07 & 0.0131 & 2.81 $10^{-3}$\\
0.05 & 0.00398 & 9.87 $10^{-4}$\\
0.04 & 0.00143 & 4.06 $10^{-4}$\\
0.03 & 0.000234 & 6.17 $10^{-5}$\\
\hline
\end{tabular}
\label{tab:3}
\caption{Dependence of the transition probability with respect to $\epsilon$.}
\end{table}

In the following Figure \ref{fig:decExp} is plotted the logarithm of $\mathbb{P}(\{ \tau_B < \tau _A\})$
with respect to $1/\epsilon$ for the values indexed  in the above table.
We obtain a straight line showing the exponential decrease of $\log(\mathbb{P}(\{ \tau_B < \tau _A\})) $ with respect to $1/\epsilon$.
\begin{figure}[h]
   \centerline{\includegraphics[width=7cm, height =6.5cm]{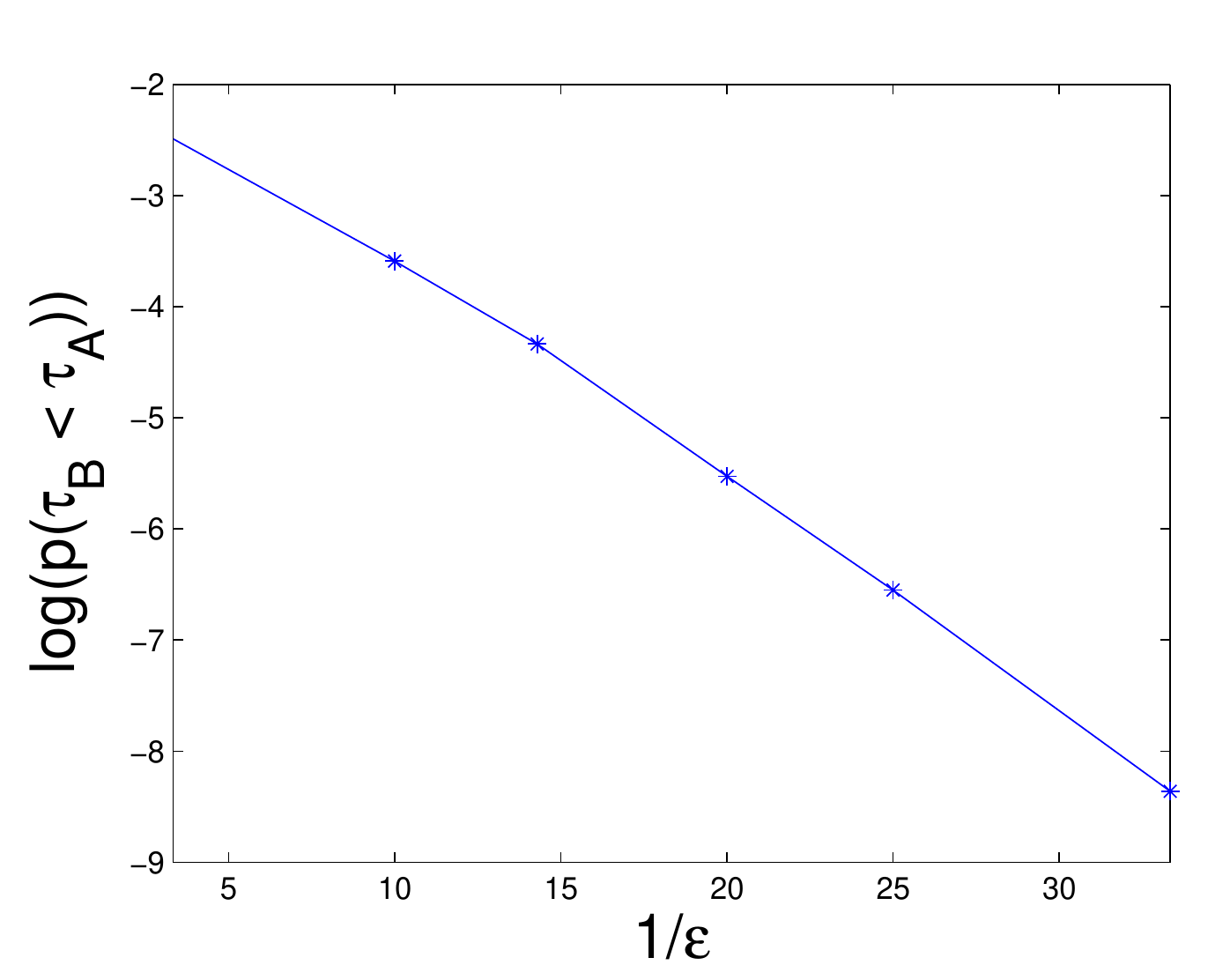}} 
	\caption{Plot of  $ \log(\mathbb{P}(\{ \tau_B < \tau _A\})) $ with respect to $1/\epsilon$.}
	\label{fig:decExp}
  \end{figure}

\subsection{Reactive trajectories in dimension $2$}
In this paragraph, we investigate numerically the qualitative behavior of the reactive trajectories in the AMS algorithm 
in terms of the parameter $\gamma$. As explained in the Section \ref{Bif}, the potential changes with $\gamma$ and then the reactive trajectories may experiment different paths to go from $A$ to $B$. At the end of the algorithm we obtain $n_{rep}$ trajectories of the process starting at the same position, and with the property that $B$ is reached before $A$.

\begin{figure}[h]
 \begin{tabular}{p{0.43\textwidth}p{0.57\textwidth}}
   \centerline{\includegraphics[width=6cm, height =5.5cm]{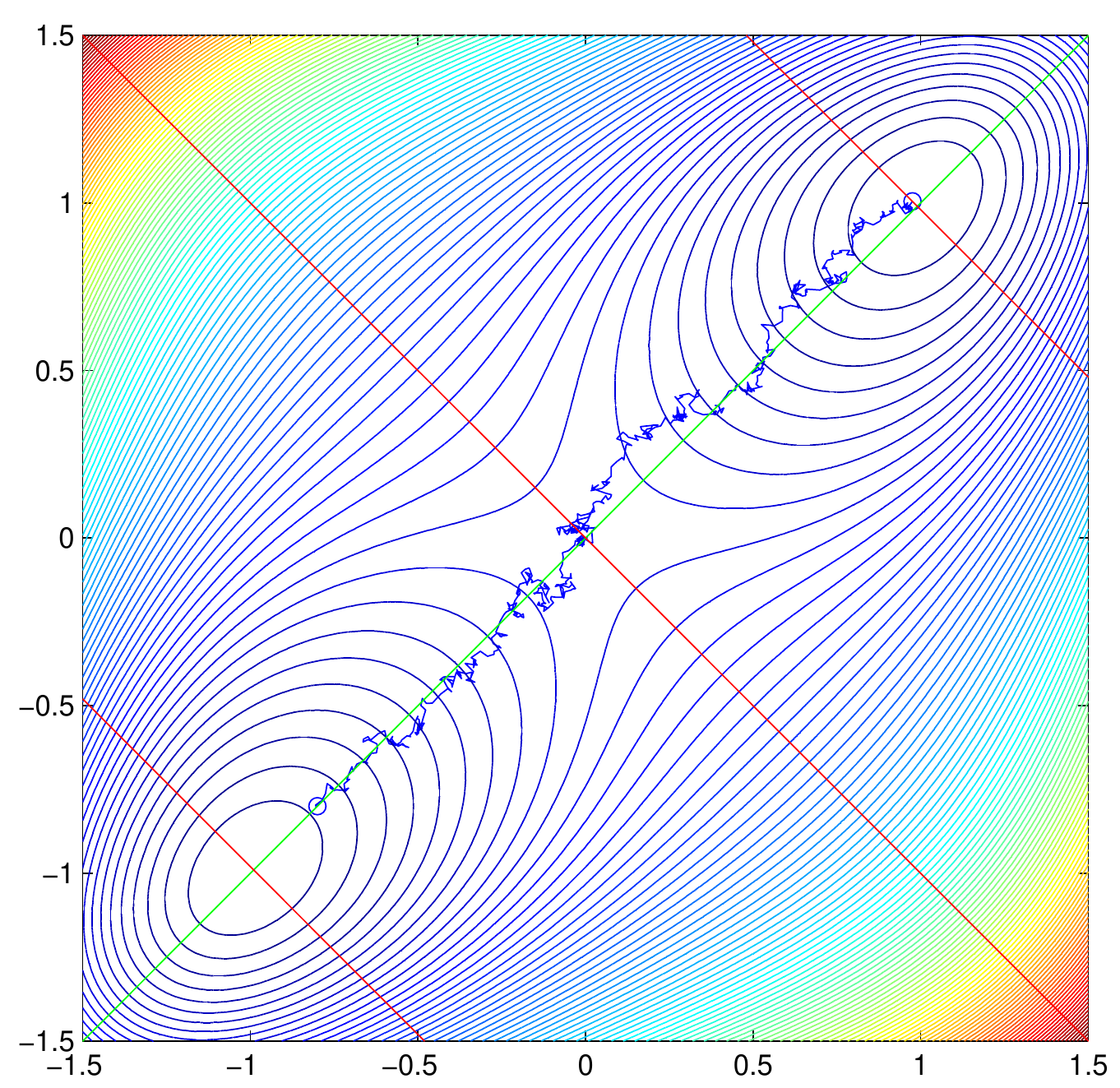}} &
   \centerline{\includegraphics[width=6cm, height =5.5cm]{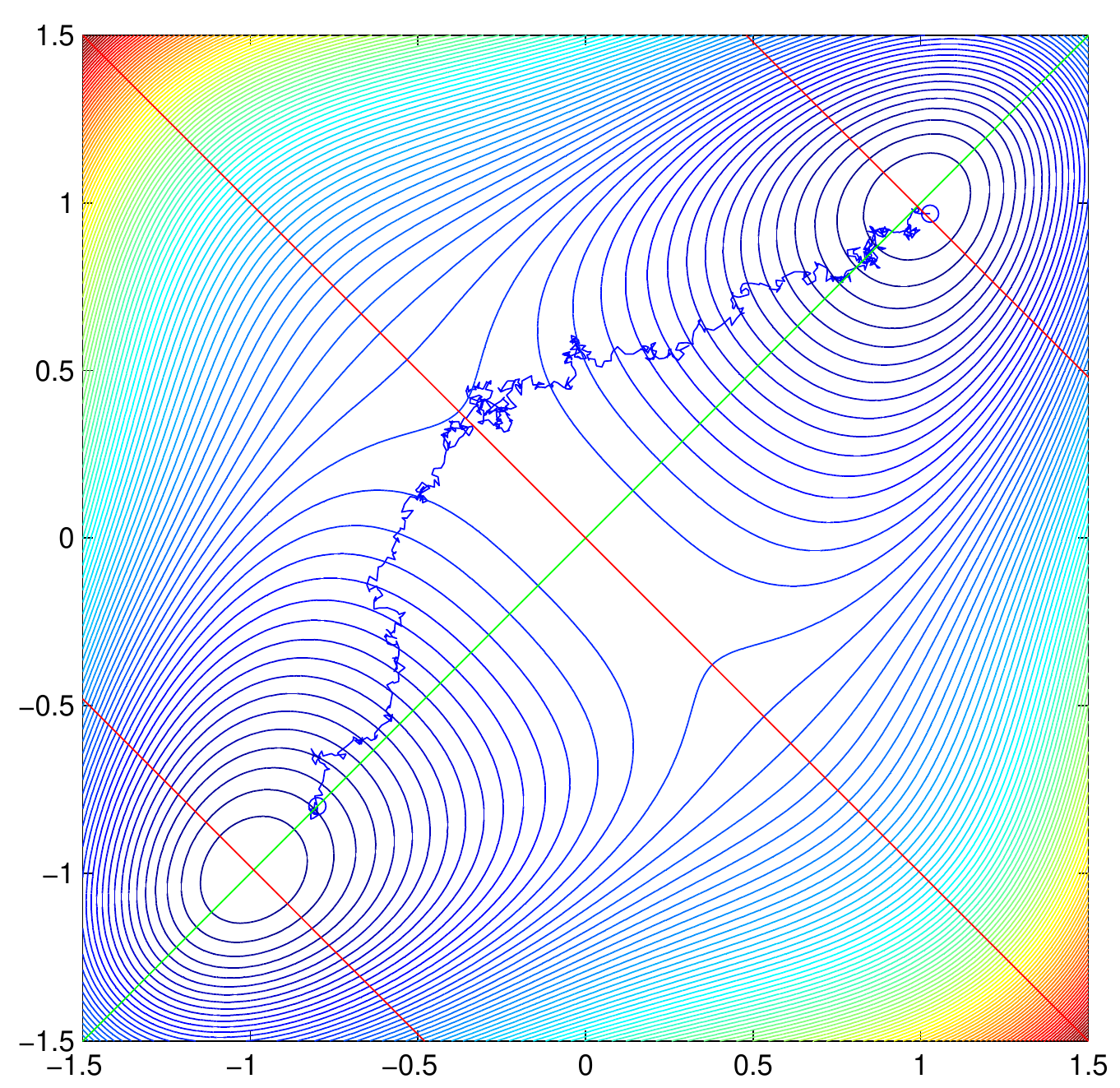}}  \\
   \centerline{\includegraphics[width=6cm, height =5.5cm]{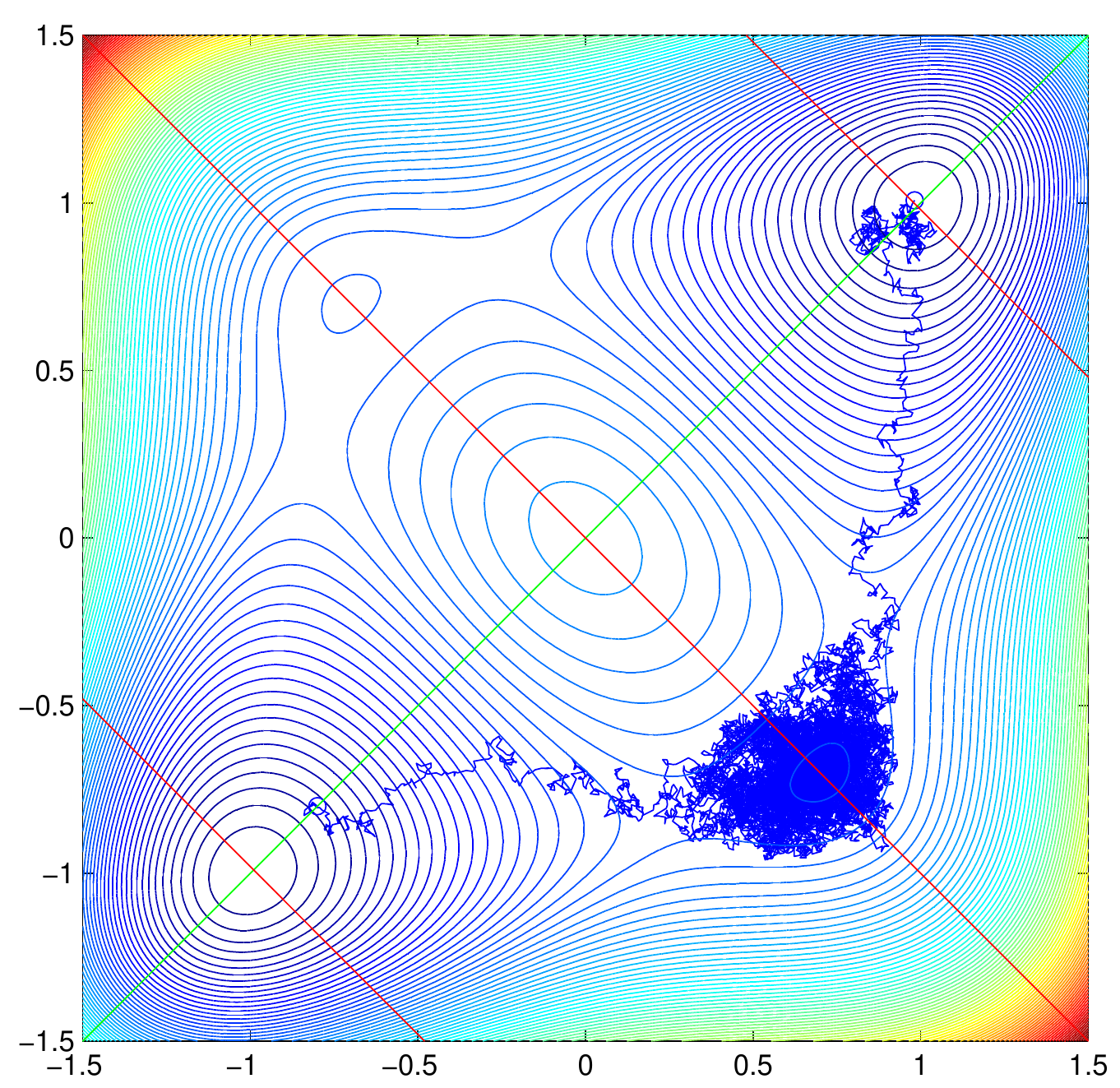}} &
   \centerline{\includegraphics[width=6cm, height =5.5cm]{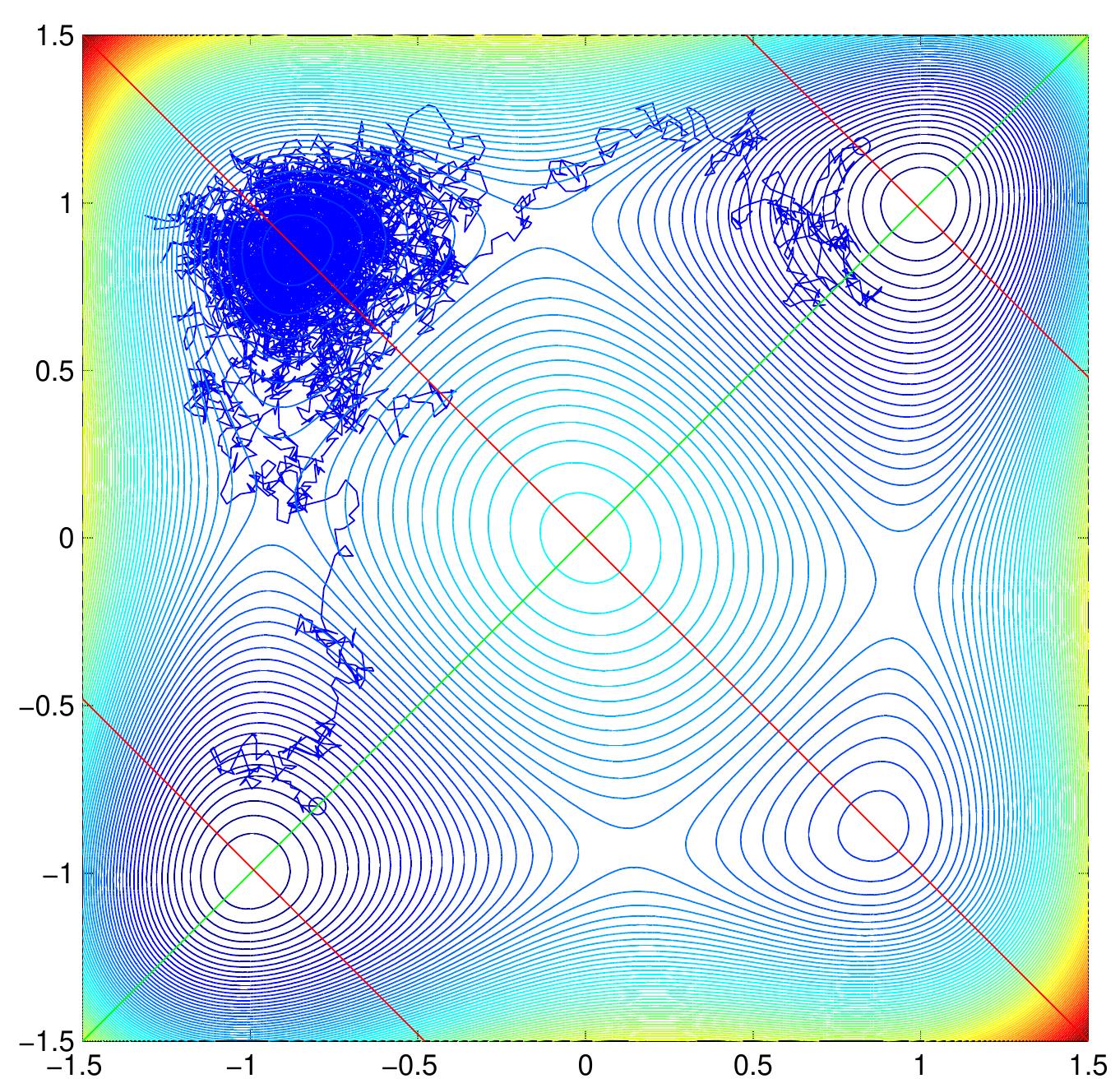}} 
  \end{tabular}
  \caption{Reactive trajectories for the two particles system and for different values of $\gamma = 1/4,1/8,1/16,1/32 $.}
\label{bif2D}
  \end{figure}

\begin{figure}[h]
 \begin{tabular}{p{0.43\textwidth}p{0.57\textwidth}}
   \centerline{\includegraphics[width=6cm, height =5.5cm]{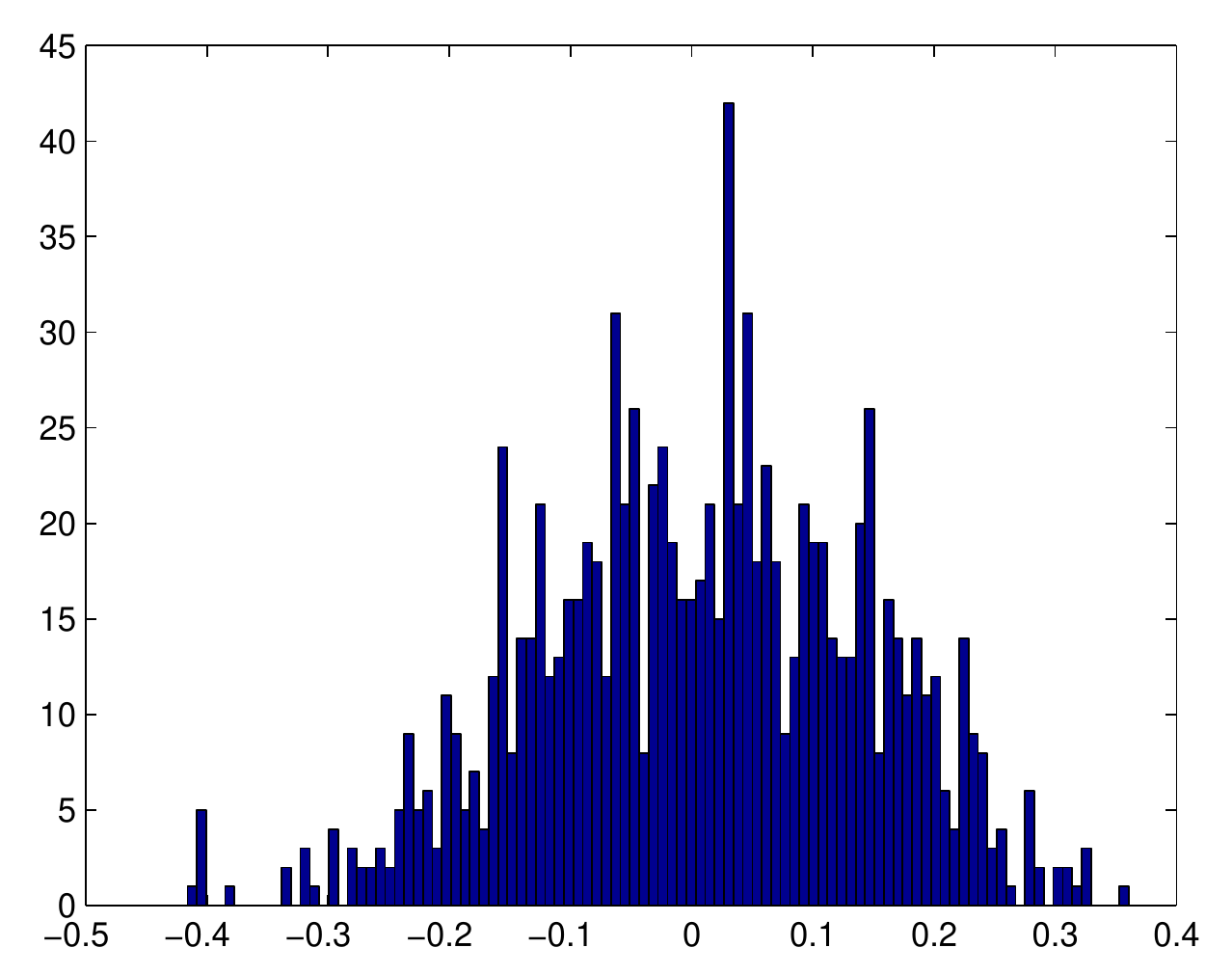}} &
   \centerline{\includegraphics[width=6cm, height =5.5cm]{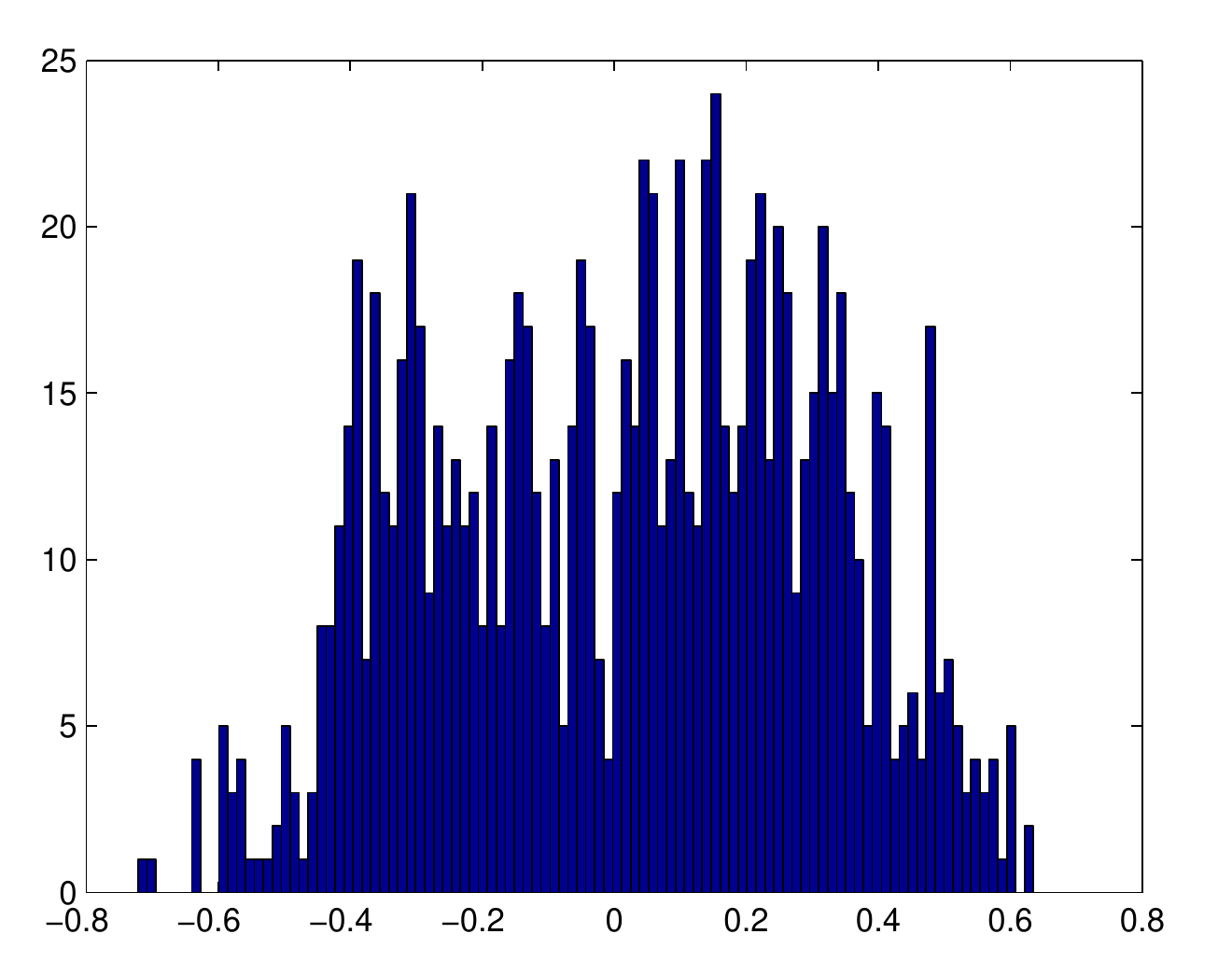}}  \\
   \centerline{\includegraphics[width=6cm, height =5.5cm]{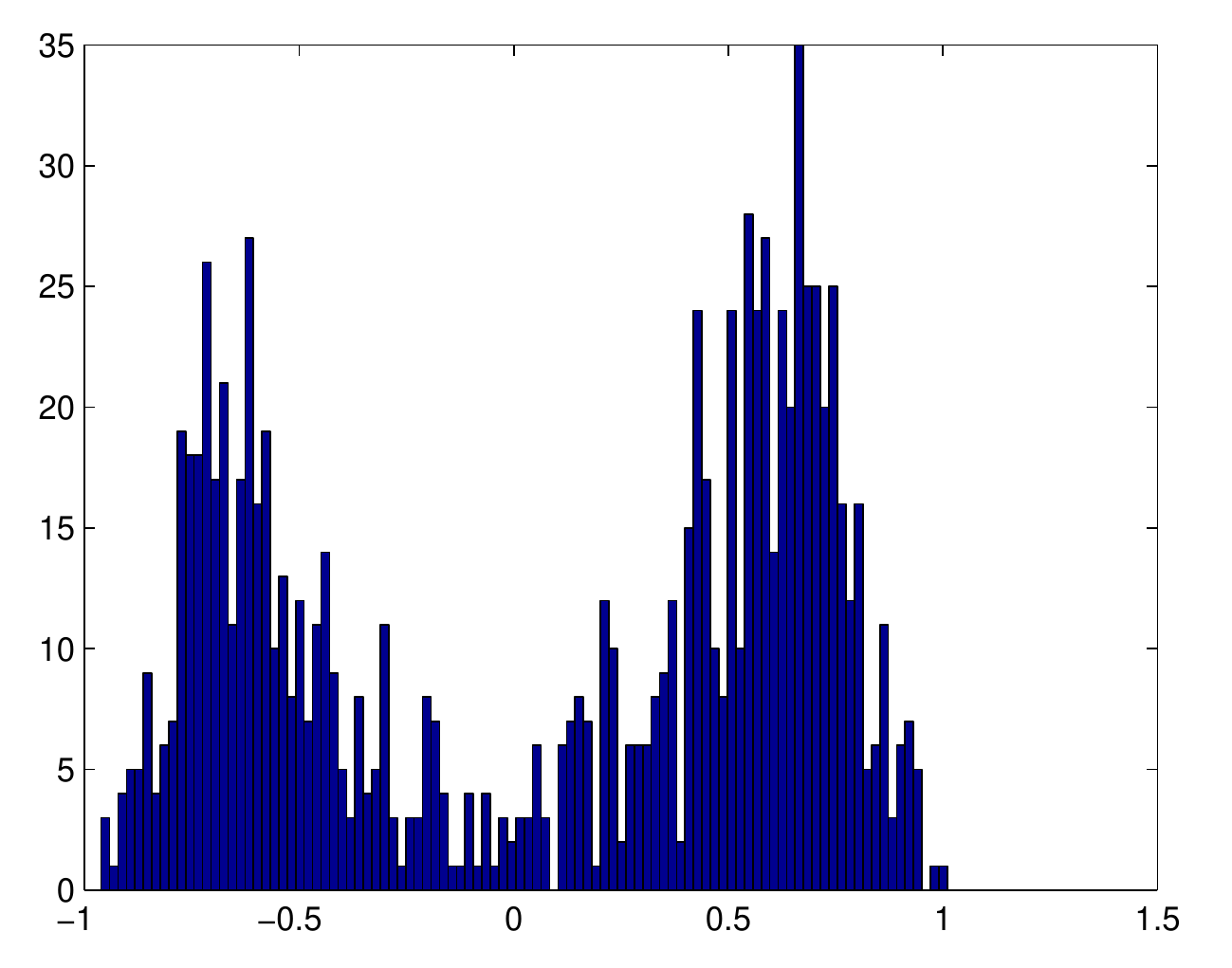}} &
   \centerline{\includegraphics[width=6cm, height =5.5cm]{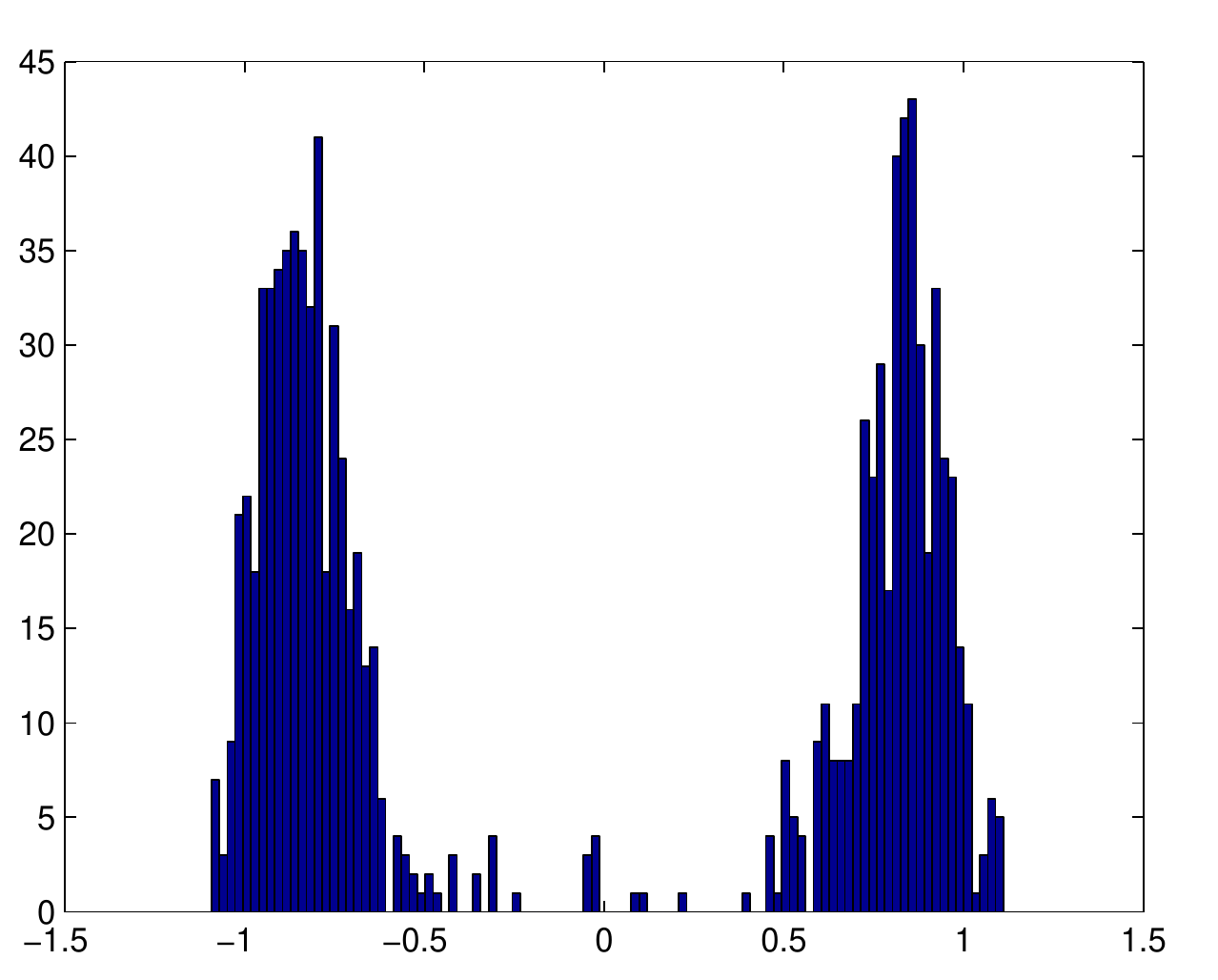}} 
  \end{tabular}
  \caption{Histograms for the two particles system in the $(x, -x)$ line and for different values of $\gamma = 1/4,1/8,1/16,1/32 $.}
\label{Hist2D}
  \end{figure}

We consider different values for $\gamma\in\left\{1/4,1/8,1/16,1/32\right\}$. For each, we have used the algorithm for two different values of the number of replica $n_{rep}$ and of the temperature $\epsilon$: either $(n_{rep}=100,\epsilon=1.10^{-4})$ or $(n_{rep}=1000,\epsilon=5.10^{-4})$. Computations take a few minutes on a personal computer. 
Figure \ref{bif2D} displays examples of reactive trajectories in dimension $N=2$ obtained with $(n_{rep}=100,\epsilon=1.10^{-4})$.
	To represent the qualitative behavior of the approximate reactive trajectories, we also draw histograms of the position in the line $x+y=0$ when reactive trajectories are crossing it. With $n_{rep}=1000$ we obtain almost symmetric histograms of Figure \ref{Hist2D}.

We recover the expected behavior with respect to $\gamma$: first, $(0,0)$ is the unique saddle point; then two better saddle points appear on each side of the line $x+y=0$; finally two local minima become places where trajectories are trapped during a long time - the importance of this effect on trajectories should decrease when temperature decreases, with the price of more iterations in the algorithm.

Notice that since we only kill $k_{rep}=1$ replica at each iteration of the algorithm, the time to either reach $A$ or $B$ can be very long even after several steps, due to the presence of additional local minima for instance. This observation is the origin for investigating parallelization strategies, to reduce computational time. Taking $k_{rep}>1$ is one natural answer to this problem, but possibly not the only one.

\subsection{Reactive trajectories for the Allen-Cahn equations}

\begin{figure}[h]
 \begin{tabular}{p{0.43\textwidth}p{0.57\textwidth}}
   \centerline{\includegraphics[width=6cm, height =5.5cm]{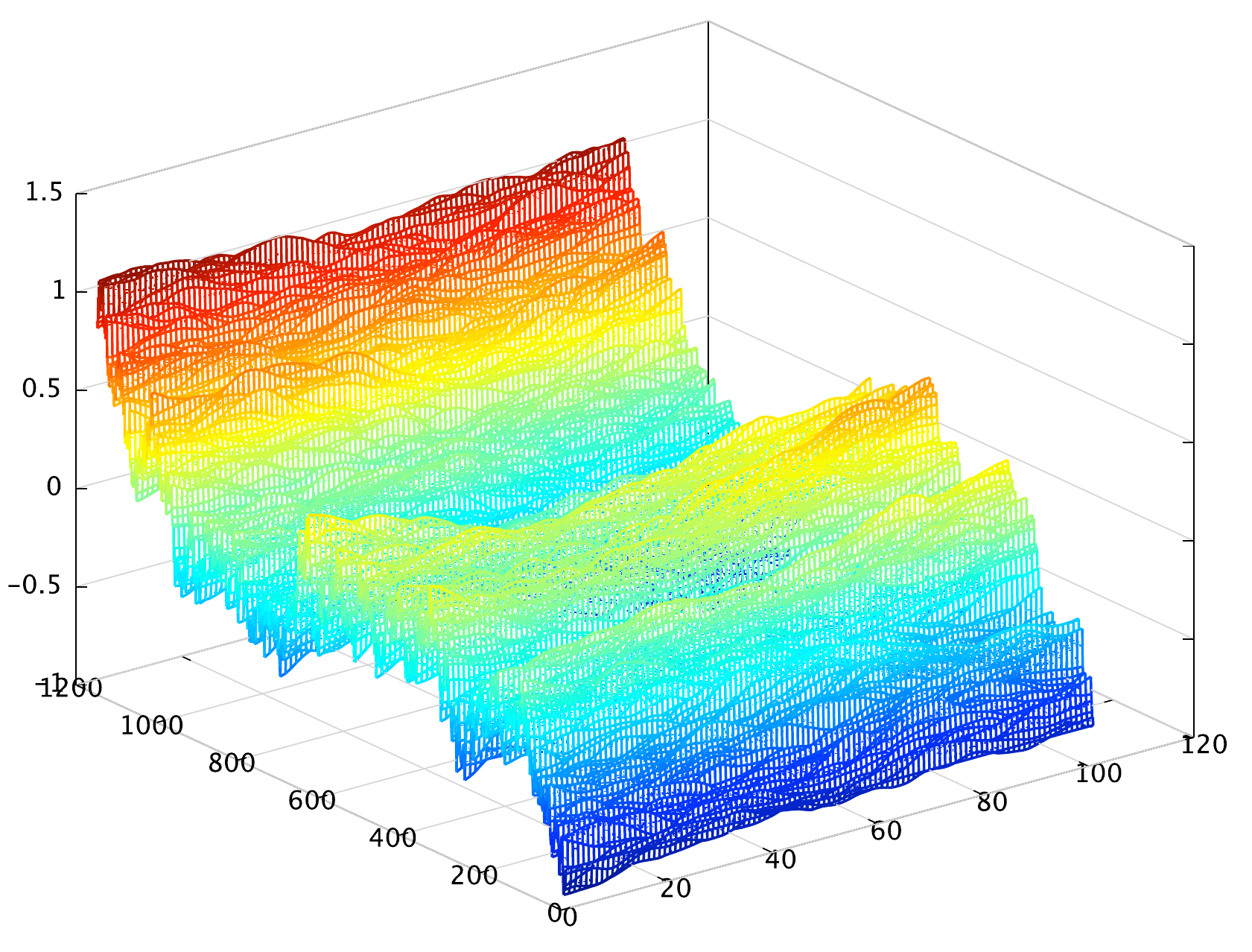}} &
   \centerline{\includegraphics[width=6cm, height =5.5cm]{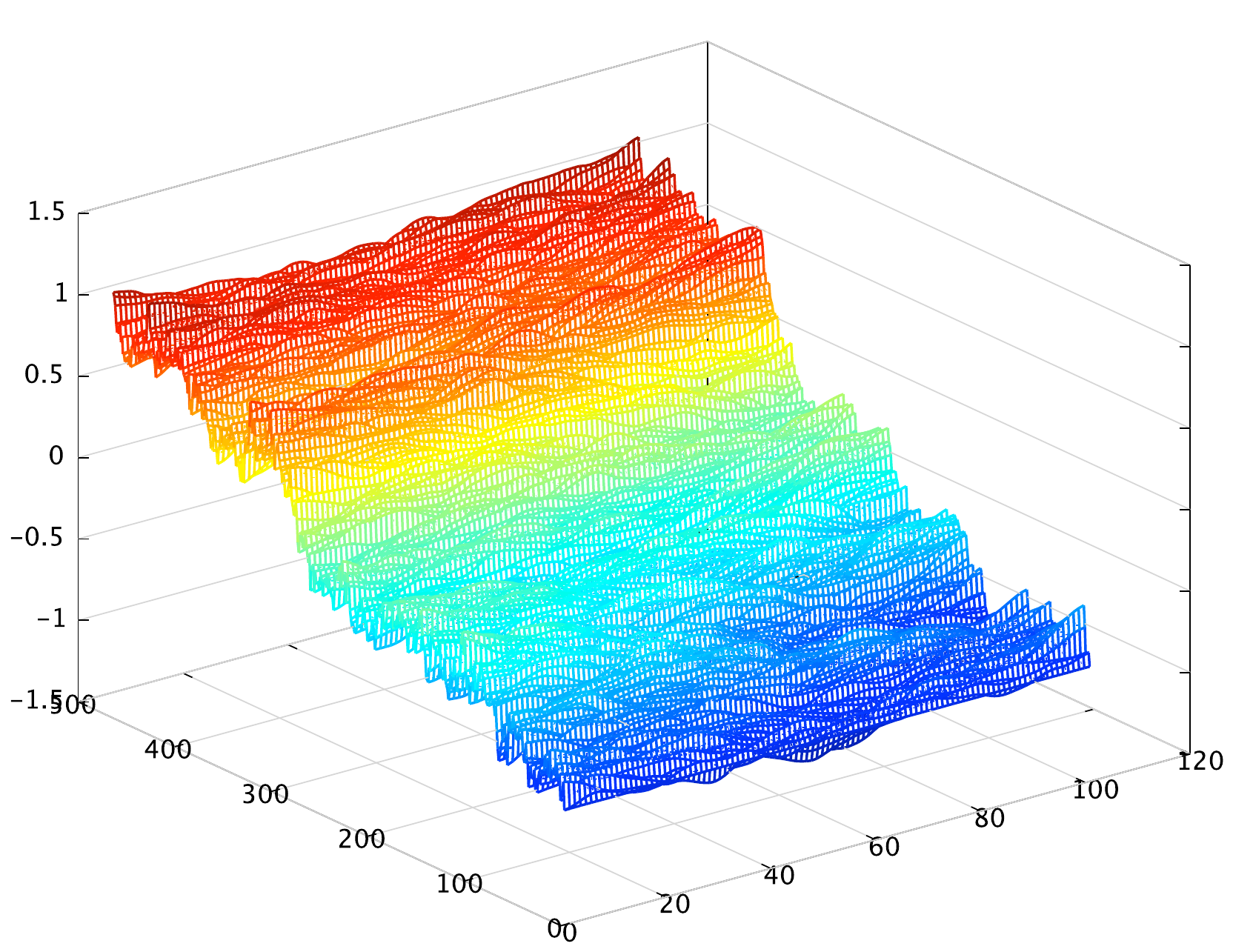}}  \\
   \centerline{\includegraphics[width=6cm, height =5.5cm]{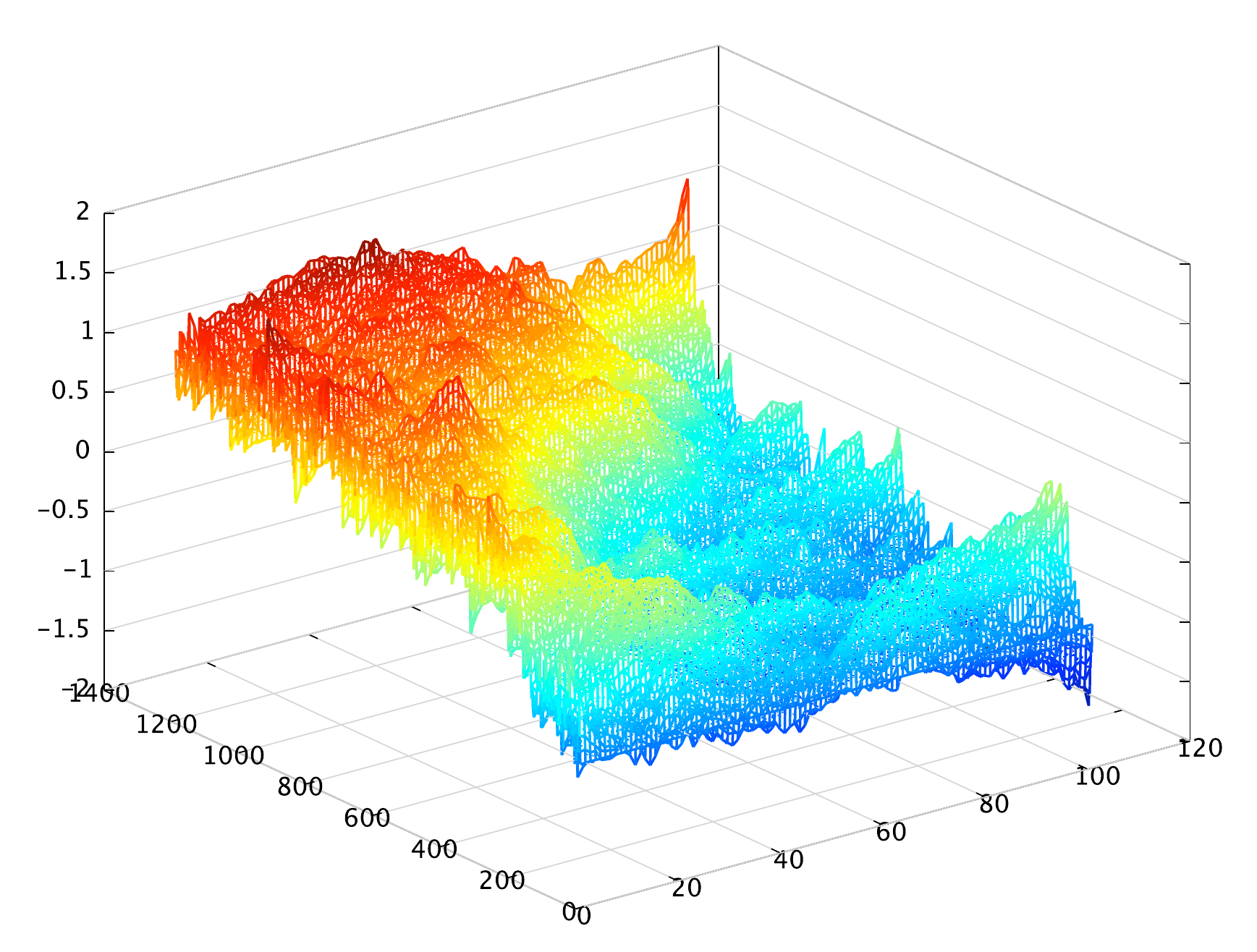}} &
   \centerline{\includegraphics[width=6cm, height =5.5cm]{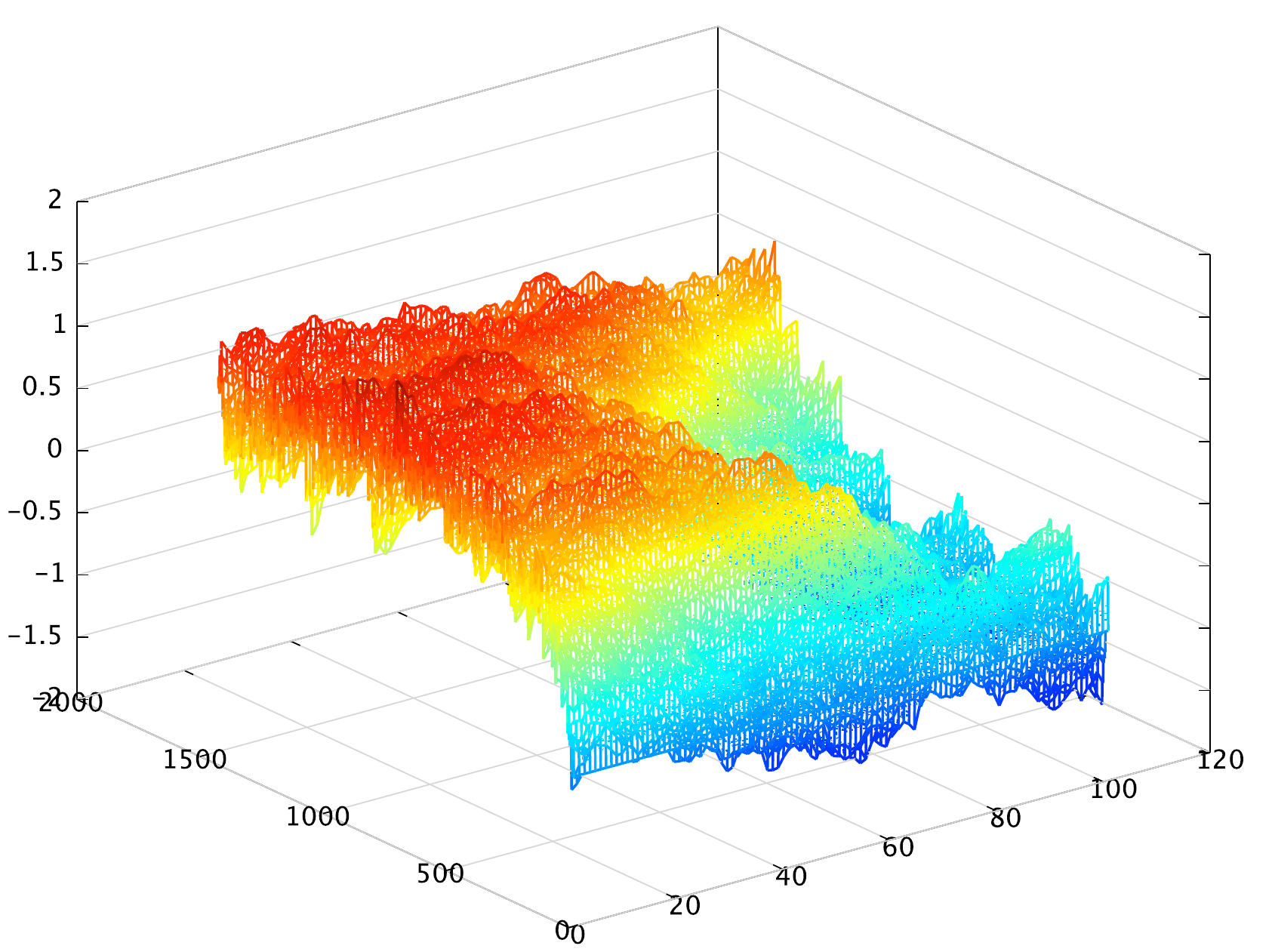}} 
  \end{tabular}
  \caption{Reactive trajectories for the Allen-Cahn equation for $\gamma=1$ and $\gamma=0.1$.}
\label{bifAC}
  \end{figure}

We give a few examples of trajectories obtained in the AMS algorithm for the Allen-Cahn equation. It is discretized as before with $\Delta t=0.01$ and $\Delta x=0.02$; we use $n_{rep}=100$ and $\epsilon=0.05$. The initial conditions is given by $x(\lambda)=-0.8$. In Figure \ref{bifAC} is given some examples of reactive trajectories for different values of $\gamma=1$ and $\gamma=0.1$.

\section{Conclusion}
In this proceeding, we investigated the numerical estimations of rare events in infinite dimensions.
We proposed an efficient algorithm outperforming the standard Monte Carlo simulations and consisting on a generalization
of splitting algorithms in finite dimension. Beyond the fact that very small probabilities can be reached, 
reactive trajectories can also be computed. These paths are known for \eqref{SystSDE} and even \eqref{SPDE1}, but it can be helpful in more complex situations. Finally, the great advantage of 
this algorithm is to be highly parallelizable. However, a lot of work remains to be done. In a forthcoming paper, 
we will study the unbiased property of the estimator for $k_{rep}\geqslant 1$.
Based on this analysis, we wish to develop a parallel version of this algorithm whose efficiency 
will depend on a ratio between the number of processors and the number of killed replicas $k_{rep}$.

\begin{acknowledgement}
The authors would like to thank the organizers of the CEMRACS 2013 (N. Champagnat, T. Leli\`evre and A. Nouy) for a very friendly research environment. The authors also acknowledge very fruitful discussions with T. Leli\`evre and D. Aristoff.
Finally, they thank D. Iampietro for its participation to the project.
\end{acknowledgement}

\bibliographystyle{plain}
\bibliography{bibAMS}
\end{document}